\documentclass{conm-p-l}
\usepackage{amssymb}

\theoremstyle{definition}

\theoremstyle{remark}

\numberwithin{equation}{section}
\newcommand{\nt}{\noindent}
\newcommand{\bs}{\bigskip}
\newcommand{\ms}{\medskip}
\newcommand{\mk}{\medskip}
\newcommand{\sk}{\smallskip}

\begin{document}

\title{Bass's Work in Ring Theory and Projective Modules}

%    Information for first author
\author{T.~Y.~Lam}
%    Address of record for the research reported here
\address{Department of Mathematics, University of California, 
Berkeley, CA 94720}
%    Current address
%\curraddr{}
\email{lam@math.berkeley.edu}
%    \thanks will become a 1st page footnote.
\thanks{The work on this paper was supported in part by a grant from NSA}

%    General info
\date{March 23, 1999.}

\dedicatory{This paper is dedicated to Hyman Bass on his 65th birthday.}

\begin{abstract}
The early papers of Hyman Bass in the late 50s and the early 60s leading 
up to his pioneering work in algebraic $\,K$-theory have played an
important and very special role in ring theory and the theory of projective 
(and injective) modules.  In this article, we give a general survey of 
Bass's fundamental contributions in this early period of his work, and 
explain how much this work has influenced and shaped the thinking of 
subsequent researchers in the area.
\end{abstract}

\maketitle

\bigskip

\medskip
\noindent {\bf Contents}

\bs
\nt {\bf \S0.  \ Introduction} \\

\noindent {\bf Part I: Projective (and Torsionfree) Modules}  \\

\nt {\bf \S1.  \ Big Projectives}  \\
\nt {\bf \S2.  \ Stable Structure of Projective Modules}  \\
\nt {\bf \S3.  \ Work Related to Serre's Conjecture} \\
\nt {\bf \S4.  \ Rings with Binary Generated Ideals: Bass Rings} \\

\noindent {\bf Part II: Ring Theory} \\

\nt {\bf \S5.  \ Semiperfect Rings as Generalizations of Semiprimary Rings}\\
\nt {\bf \S6.  \ Perfect Rings and Restricted DCC}    \\
\nt {\bf \S7.  \ Perfect Rings and Representation Theory}  \\
\nt {\bf \S8.  \ Stable Range of Rings}  \\
\nt {\bf \S9.  \ Rings of Stable Range One}

\bs
\nt {\bf References}

\bs

\mk
\begin{center}
{\bf \S0. \  Introduction}
\end{center}
\medskip

It gives me great pleasure to have this opportunity to write about
Professor Hyman Bass's work in ring theory and projective modules.  This was
work done by a young Hyman in the early 60s when he was a junior faculty
member at Columbia.  The now-classical paper on the homological generalization
of semiprimary rings [B$_1$], an outgrowth of his 1959 Chicago thesis written 
under the direction of Irving Kaplansky, was followed in quick succession 
by a brilliant series of papers [B$_2$--B$_7$] on various aspects of the 
structure of projective modules (decomposition, extendibility, and freeness) 
and injective modules (injective dimensions, Gorenstein rings, and Bass 
rings).  The 1962 announcement of a homotopy theory of projective modules 
with Schanuel [BS] and the subsequent 1964 announcement [B$_8$] on the 
stable structure of the general linear group over an arbitrary ring 
culminated in his famous IHES paper [B$_9$] which, perhaps more than any 
other work in that era, marked the monumental creation of the new 
mathematical subject of algebraic $\,K$-theory.  The usual way to put it 
would be to say, I guess, that all of this work ``showed the future master''; 
but in fact, by this time, Hyman Bass had already proved himself to be a 
true master of the great art of algebra --- and a mathematician of 
extraordinary creativity and insight.

\bs
Since Chuck Weibel [We] will be touching upon Bass's work in algebraic 
$\,K$-theory and Craig Huneke [Hu] will be reporting on his work on injective 
modules and Gorenstein rings, I shall focus the discussion in this article 
mainly on Bass's contributions to ring theory and the structure of projective 
modules. Starting with the latter, we first survey, in \S1 and \S2 below, 
Bass's work on big projectives and on his generalizations and extensions of 
Serre's results on the stable structure of projective modules.  In \S3, we 
move on to Bass's work surrounding the theme of Serre's Conjecture on 
f.g.\footnote{Throughout this paper, ``f.g.'' will be used as an abbreviation
for ``finitely generated''.} projective modules over polynomial rings. The 
last section of Part I concludes with a survey on Bass's work on the
decomposition of torsionfree modules and on commutative rings with binary 
generated ideals, which resulted in the important notion of Bass rings 
and Bass orders.  In this section (\S4), Gorenstein rings play a
substantial role, but our exposition is designed to overlap only 
peripherally with that of Huneke's article.

\bs
In Part II of this article, we return to Bass's maiden work on the homological 
generalizations of semiprimary rings, and explain the historic context of this
work as well as the role it played in the subsequent development of 
noncommutative ring theory.  Guided by the ideas of stability in the homotopy 
theory of vector bundles, Bass also single-handedly invented the notion of 
stable range of a ring $R$, which he successfully applied to the study of 
the stable structure of the infinite general linear group over $R$.  (This 
work may be thought of as the ``$K_1\,$-analogue'' of the stable theory of 
projective modules reported in \S2.)  Our survey concludes with a discussion 
(\S\S8-9) of Bass's work on the stable range, with a special emphasis on the 
case of stable range one which in turn has deep significance on the 
arithmetic of rings, and on questions concerning the cancellation and 
substitution of modules with respect to direct sum decompositions.

\bs
Throughout this exposition, our aim is not only to survey Bass's work 
in ring theory and projective modules, but also to point out how much
this work has influenced and shaped the thinking of subsequent researchers 
in the area.  It will be seen that, in a number of new lines of
investigation in algebra in the last few decades, it was the decisive
pioneering steps of Bass that broke open the new path. Although Bass's 
work in modules and rings spanned only the decade of the 60s, its impact 
has been enormous indeed, and will certainly continue to be felt as ring 
theory moves into the next century.

\bs
Writing this article is for me a very pleasant task, though undertaking 
such a task inevitably involves somewhat of a nostalgic trip down the 
memory lane.  Since it is perhaps not out of place in these Proceedings 
to talk about one's connections to Hyman, I will indulge myself in a few 
personal reminiscences below.  

\bs
As a beginning graduate student at Columbia in the mid-60s, I was more 
than a little awed by the power and fame of the senior faculty: 
Professors Lipman Bers, Samuel Eilenberg, Ellis Kolchin, Masatake
Kuranishi, Serge Lang, Edgar Lorch, and Paul A.~Smith, among others, from 
whom I took my first and second year graduate courses.\footnote{The students 
were nothing short of stellar either!  My graduate peers included Mike 
Engber, Audun Holme, Fred Gardiner, and Irwin Kra; Winfried Scharlau was a 
visiting student from Germany, and Alexander Mikhalev was an exchange scholar
from U.S.S.R.; Bill Haboush, Tony Bak, Mike Stein, Allen Altman, Bruce Bennett,
Spencer Bloch, Bob and Jane Gilman were a couple of years behind; Julius 
Shaneson, Sylvain Cappell and Ethan Akin were ``hot-shot'' Columbia 
undergraduates taking graduate courses with us. Each person in this list is 
now a Professor of Mathematics.  Many of them have been chairs, provosts, and 
deans.}  After a few skirmishes with functional analysis (taught to me by 
Professor Lorch and Visiting Professor B.~Sz.~Nagy), I fell under the 
spell of Sammy Eilenberg, and positioned myself to become a student of his, 
hoping to study with him category theory and homological algebra.  Under 
Sammy's guidance, I wrote my first paper [La$_1$], which he kindly 
communicated to the Proceedings of the National Academy of Sciences (of 
which he was a member). But then for the year 1966-67, Professor Eilenberg 
had his Sabbatical coming, and he was to go off to Paris to spend the entire 
academic year.  To go with him to Paris would have been the ``cool'' 
thing to do for a mid-career graduate student except, alas, I found myself 
both financially and linguistically poorly equipped to make the trip.  So 
Professor Eilenberg said to me: ``Why don't you stay here, and study with 
Bass while I am gone?''  This was how I became a student of Bass!  Hyman was 
then a young Assistant Professor, who, just a few years ago, was brought 
from Chicago to Columbia by Sammy himself.

\bs
Hyman became an Associate Professor around that time, and in another 
year he was moved up to Full Professor.  Even I could see he must
be working on something hot!  So I xeroxed all of his papers,\footnote{Yes,
there were already xerox machines in the mid-60s, though stencils and 
mimeographed copies were still not entirely out of fashion.} and started 
poring over them assiduously.  Hyman's $\,K$-theory paper [B$_9$] had just 
come out in the ``Blue Journal''; this and his earlier papers in ring theory 
and projective modules eventually became a staple of my graduate education.  
The most challenging open problem in homological algebra in those days was 
Serre's Conjecture ({\it ibid.}); many a graduate student in algebra from 
that era had no doubt tried his/her hand at it --- and I was no exception.  
Hyman, who had generalized Seshadri's solution of the Conjecture in two
variables, was the natural leader for this small circle of aspiring
graduate students. I still recall that, one time, one of us thought 
he had a brilliant idea to solve Serre's Conjecture.  A few of us 
excitedly met with Hyman in an impromptu seminar to go over the ``idea''; 
but of course Hyman quickly found the hole.  Although we never got 
anywhere with our fledgling efforts, my fascination with Serre's 
Conjecture continued, and culminated in the writing of my 1978 Springer 
Lecture Notes [La$_2$], two years after the Conjecture was fully solved 
independently by Suslin and Quillen.  Hyman's great influence was evident 
from cover to cover of my modest book.

\bs
I spent my last graduate year at Columbia in 1966-67.  In that year,
Hyman taught what was most probably the first graduate course ever given
in algebraic $\,K$-theory in the US.  Pavaman Murthy had come from the 
Tata Institute to do postdoctoral work with Hyman, and was in the audience.
(In his modest Morningside Heights apartment, Pavaman served what he 
claimed to be the best coffee in Manhattan.  It was free, so I had no 
reason to disagree; we became good friends.)  Some of the lecture 
notes taken by Murthy, Charles Small and me eventually evolved into Hyman's 
famed tome [Ba$_{10}$].  At that time, Hyman was busy at work with Milnor 
and Serre on the Congruence Subgroup Problem for the special linear groups 
and the symplectic groups.  I got lucky and proved several little things about 
Mennicke symbols and $\,SK_1\,$ of abelian group rings, which earned me a 
few attributions in [BMS].  Needless to say, I was proud to be mentioned 
in a paper by such distinguished authors.  I got lucky in some other fronts 
too, and in May, 1967 completed a thesis in algebraic $\,K$-theory under 
Hyman dealing with Artin's Induction Theorem and induction techniques for 
Grothendieck groups and Whitehead groups of finite groups.  In those days, 
algebraic $\,K$-theory meant only $\,K_0\,$ and $\,K_1$; even Milnor's 
$\,K_2\,$ had not been defined yet.

\bs
As it happened, I was Hyman's first Ph.D.~student.  I have always viewed
this as a special honor, and I am sure that this fact has helped me a 
great deal professionally.  Now the list of Ph.D.~students of Hyman is 25(?)
strong.  The longevity of Hyman as a thesis advisor is rather strikingly
illustrated by the fact that at least several of my mathematical brothers 
and sisters in this list were not even born yet when I completed my 
Ph.D.~degree at Columbia.

\bs
In supervising my work, Hyman never tried to tell me what I should work
on.  Rather, he let me find my own path, and instilled in me the needed 
confidence to grow into a research mathematician.  What I learned from
Hyman was not just mathematics, but how to do mathematics, and, what
is perhaps even more important, how to conduct myself as a mathematician.
All of this he taught me in the best way --- by his own example.  For 
this, and for the many other favors he has rendered over the years, 
I shall always be grateful.

\mk

\bs
\begin{center}
{\bf Part I: Projective (and Torsionfree) Modules}  

\bs
{\bf \S1. Big Projectives}
\end{center}

\sk
While the projective modules occurring in number theory, representation
theory and algebraic geometry are mostly f.g.~ones, non-f.g.~projective 
modules do arise naturally over various kinds of rings, for instance, 
rings of continuous functions.  Thus, the quest for information about 
non-f.g.~projective modules is not a frivolous one.

\bs
The first significant result in the study of general (that is, not
necessarily f.g.) projective modules was found by Kaplansky
in 1958.  In his seminal paper [K$_2$], Kaplansky proved that, for any ring
$R$, any projective $R$-module is always a direct sum of countably
generated (projective) modules.  This result has two interesting consequences.
First, any indecomposable projective module over any ring is countably 
generated; and second, if any countably generated projective module over
some ring $R$ is free, then {\it any\/} projective $R$-module is also free.
The case of a local ring $R$ provides a particularly striking illustration 
for the power of the second statement: in this case, Kaplansky used an 
ingenious argument (reminiscent of the proof of Nakayama's Lemma) to show 
that any countably generated projective $R$-module is free, from which 
it then follows that any projective $R$-module is free. 

\bs
Inspired by Kaplansky's result, Bass took up the study of ``big''
projective modules in [B$_5$].  The overall theme of [B$_5$] is that,
under certain fairly mild conditions on a ring $R$, ``big'' projective 
$R$-modules are necessarily free.  To formulate this precisely, Bass
introduced the notion of a {\it uniformly $\,\aleph$-big\/} module: 
for an infinite cardinal $\,\aleph$, an $\,R$-module $\,P\,$ is uniformly
$\,\aleph$-big if $P$ can be generated by $\,\aleph\,$ elements, 
and for any ideal $\,I\subsetneq R$, the $\,R/I$-module 
$\,P/I\,P\,$ cannot be generated by fewer than $\,\aleph\,$ elements.  
For instance, the free module of rank $\,\aleph\,$ over a nonzero ring $R$ is 
uniformly $\,\aleph$-big.  Bass's first result on big projectives is
the following converse of this statement:

\bs\nt
{\bf (1.1) Theorem.} {\it If} $\,R/\mbox{rad}(R)\,$ {\it is left 
noetherian, then, for any infinite cardinal $\,\aleph$, any uniformly} 
$\,\aleph$-{\it big projective left $R$-module is free.}

\bs
Bass's proof of this result is modeled upon Kaplansky's proof of his main
theorem in [K$_2$]: one makes a reduction to the crucial case when $\,\aleph\,$
is the countable infinite cardinal $\,\aleph_0$, and in this case one achieves
the desired goal by a clever juggling with infinite matrices.  

\bs
For concrete applications of (1.1), one needs to find situations where we 
can say, for instance, that all non-f.g.~projective modules 
are uniformly big for some infinite cardinal.  Bass showed that this is 
the case when, say, $R$ is a commutative noetherian ring with only trivial 
idempotents.  Thus, one has the following 

\bs\nt
{\bf (1.2) Corollary.} {\it If $R$ is a commutative noetherian ring with
only trivial idempotents, then any non-f.g.~projective $R$-module is free.}

\bs
For instance, if $\,R=k[x_1,\dots,x_n]\,$ where $k$ is a field, Serre
asked in 1955 (see \S2) if every f.g.~projective $R$-module $P$ is free. 
The result above would give the freeness of $P$ if $P$ was {\it not\/} f.g.

\bs
The general results of Bass (such as (1.1) and (1.2)) describing the behavior 
of big projectives have remained essentially unsurpassed to this date.
In [B$_5$], Bass remarked that these results seemed to indicate that 
big projective modules ``invite little interest.''  We can say, today, that
this is perhaps not quite true.  In a recent paper [LR], Levy and Robson have
determined the structure of all infinitely generated projective modules over 
(noncommutative) hereditary noetherian prime rings.    Their results showed
that, over such rings, there may exist non-f.g.~projective modules which
are not free and which have rather interesting structures.

\bs
Two other results of Bass on projective modules over any ring have become 
folklore in the subject, so it is fitting to end this section by recalling
them.  The first of these says that:

\bs\nt
{\bf (1.3) Theorem.} {\it Any nonzero projective module $P$ over any ring $R$
has a maximal submodule.  More precisely,} rad$(R)\,P=\mbox{rad}(P)\neq P$,
{\it where\/} rad$(P)$ {\it denotes the intersection of all maximal submodules
of $\,P$.}

\bs
This result first appeared in [B$_1$:~p.\,474], and has been used by many 
authors since.  The second result, also involving the Jacobson radical
of a ring $R$, states the following:

\bs\nt
{\bf (1.4) Theorem.} {\it Let $\,I\,$ be any ideal of a ring $R$ such
that\/} $I\subseteq \mbox{rad}(R)$, {\it and $\,P,\,Q\,$ be f.g.~projective 
left $\,R$-modules.  Then 
$$ P\cong Q \;\,as \,\;R\mbox{-}modules \;\Longleftrightarrow
P/IP\cong Q/IQ \,\; as \,\;R/I\mbox{-}modules. \leqno (1.5) $$
In particular, $P$ is free
as an $R$-module iff $\,P/IP\,$ is free as an $\,R/I$-module.}

\bs
This result, which made its first appearance in the literature as Lemma 2.4 
in [B$_2$], is now in virtually every textbook which treats the subject
of projective modules, and is well-known to any student knowledgeable
about the subject of homological algebra.  We should point out that, 
with Bass's notion of left $\,T$-nilpotency on ($1$-sided) ideals 
(introduced later in \S6), (1.4) can also be given the following 
``transfinite'' formulation, for {\it any\/} pair of projective modules
$\,P,\,Q$:

\bs\nt
{\bf (1.6) Theorem.} {\it  If $\,I\,$ is any left $T$-nilpotent 
left ideal in a ring $R$, then $(1.5)$ holds for any (not necessarily 
f.g.) projective left $R$-modules $\,P,\,Q$.}

\bs
For a proof of this, see [La$_3$:~(23.17)].

\bs
\sk
\begin{center}
{\bf \S2.  \ Stable Structure of Projective Modules}  \\
\end{center}

\sk
In his epoch-making papers [S$_1$, S$_2$], Serre established an
analogy between projective modules (in algebra) and vector bundles (in
topology). This important analogy, which was further promulgated by the 
work of Swan [Sw$_1$], enables one to establish the rudiments of a
dictionary to translate the language of projective modules into that
of vector bundles.  Now by the late 1950s, the topology of vector bundles 
was already rather well-developed; the Serre-Swan analogy mentioned above 
made it possible, therefore, for algebraists to ``predict'' (if not prove)
theorems in the theory of projective modules based on their knowledge of 
topological results in the theory of vector bundles.

\bs
One of the best known elementary facts about vector bundles is the 
following.  If $\,[X,Y]\,$ denotes the set of homotopy classes of maps 
from $X$ to $Y$, and BO$(r)$ denotes the classifying space of the 
orthogonal group O$(r)$, then

\bs\nt
{\bf (2.1) Theorem.} {\it For any connected finite\/} CW {\it complex $X$ 
of dimension $\,d$, the natural map}
$$  i_r: \,[X,\mbox{BO}(r)] \longrightarrow [X,\mbox{BO}(r+1)] $$ 
{\it is surjective for $\,r\geq d$, and injective (and hence bijective)
for $\,r\geq d+1$.}

\bs
Now by the classification theorem of vector bundles, $\,[X,\mbox{BO}(r)]\,$
represents the set of equivalence classes of real $\,r$-plane bundles over
$X$, and, with this interpretation, the map $\,i_r\,$ in (2.1) above is given
by adding a trivial line bundle.  To make the transfer into algebra, we 
replace $X$ by a commutative noetherian ring $\,R\,$ with only trivial 
idempotents (so that the Zariski prime spectrum Spec$(R)$ is a connected 
space), and replace $\,[X,\mbox{BO}(r)]\,$ by the set $\,{\mathcal P}_r(R)\,$ 
of isomorphism classes of f.g.~projective $R$-modules of rank $\,r$.  
The $\,i_r\,$ in (2.1) is then to be replaced by the map given by 
``adding a copy of $\,R$''.  With such a transfer in place, Serre broke 
new ground by coming up with commutative algebra techniques to prove the 
following remarkable analogue of the surjectivity part of (2.1).

\bs\nt
{\bf (2.2) Theorem.} ([S$_2$: Th\'eor\`eme 1]) {\it In the above setting, 
assume that the maximal ideal spectrum} \,max$(R)\,$ {\it (with the Zariski
topology) is a space of dimension\footnote{The {\it dimension\/} 
of a topological space $X$ is defined to be the supremum of the codimensions 
of its nonempty closed sets.  Here, the codimension of an irreducible closed 
set $F$ is the supremum of the lengths of finite chains of irreducible closed 
sets ``above'' $F$, and the codimension of an arbitrary closed set $C$ is 
the infimum of the codimensions of all irreducible closed sets $F\subseteq C$.}
$\,d$.  Then the map 
$$ i_r: \,{\mathcal P}_r(R)\rightarrow \mathcal P _{r+1}(R)    \leqno (2.3) $$ 
is surjective if $\,r\geq d$.  In other words, any f.g.
projective $\,R$-module of rank $\,p> d\,$ is of the form 
$\,R^{p-d} \oplus P_0\,$ for some (projective) module $\,P_0$.}

\bs
In this setting (and under the same hypotheses), the question whether the 
map $\,i_r\,$ in (2.3) is {\it injective\/} for $\,r\geq d+1\,$ (in analogy 
to the second part of (2.1)) begs to be asked.  This amounts to the following 
cancellation question for f.g.~projective $\,R$-modules 
$\,P,\,Q\,$ and $\,M$: {\it if $\,M\oplus P \cong M\oplus Q\,$ and 
$\,P\,$ has rank $>d$, does it follow that $\,P\cong Q$?}  
It would be difficult to imagine that Serre was not aware of this very 
natural question in 1957-58, but anyhow, Serre did not pursue it in [Se$_1$].

\bs
In his two papers [B$_2$] and [B$_9$] (see also the research announcement
[BS] written jointly with S.~Schanuel), Bass not only answered this question 
in the affirmative, but also relaxed some of the assumptions in (2.2), and 
extended all results to a {\it noncommutative\/} setting. Working now with a 
commutative ring $A$ and a module-finite $A$-algebra $R$, Bass considered 
$\,R$-modules $P$ which need not be projective or f.g.  He defined such a 
module to be of $\,f$-rank $\geq r\,$ if, at every maximal ideal 
$\,\mathfrak m\in \mbox{max}(A)$, the $\,R_{\mathfrak m}$-module 
$\,P_{\mathfrak m}\,$ contains an $\,R_{\mathfrak m}$-free direct summand of 
rank $\,r$.

\bs\nt
{\bf (2.4) Theorem.} {\it In the above (noncommutative) setting, assume that\/}
max$(A)\,$ {\it is a noetherian space of dimension $\,d$.}

\sk\nt
(1) {\it (Splitting) If an $R$-module $P$ is a direct summand of a direct 
sum of finitely presented $R$-modules and has f-rank $>d$, then $\,P\,$ 
has a direct summand isomorphic to $R$.}

\nt
(2) {\it (Cancellation) Let $\,P,\,Q\,$ be $R$-modules such that $P$ has a 
projective direct summand of f-rank $> d$.  Then, for any f.g.~projective
$R$-module $M$,}
$$ M\oplus P \cong M\oplus Q \Longrightarrow P\cong Q. \leqno (2.5) $$

\bs
If we now define $\,{\mathcal P}_r(R)\,$ to be the set of isomorphism classes
of f.g.~projective $R$-modules of $\,f$-rank $r$, clearly
(1) and (2) above imply that:

\bs\nt
{\bf (2.6) Corollary.} {\it The map $\, i_r: {\mathcal P}_r(R)\rightarrow 
{\mathcal P}_{r+1}(R)\,$ (defined by adding $R$) is surjective for 
$\,r\geq d$, and injective (and hence bijective) for $\,r\geq d+1$.}

\bs
Just as (2.1) is quantitatively the best result for bundles in topology, 
both (2.2) and (2.5), (2.6) are quantitatively the best for f.g.~projective
modules in algebra.  We'll mention the usual examples to substantiate this 
statement in the commutative and noncommutative cases below.

\bs
In the commutative case, let $\,R\,$ be the real coordinate ring of the 
sphere $\,S^2$, so $\,R={\mathbb R}[x,y,z]$, with the relation
$\,x^2+y^2+z^2=1$.  This is a noetherian domain of (Krull) dimension 2.
Let $\,P=\mbox{ker}(\varphi)\,$ where $\,\varphi\,$ is the 
$\,R$-epimorphism $\,R^3\longrightarrow R\,$ given by $\,\varphi(e_1)=x,\;
\varphi(e_2)=y$, and $\,\varphi(e_3)=z$.  We have $\,R^3\cong R\oplus P$, so 
$\,P\,$ is a f.g.~projective $\,R$-module of rank 2. Here, $\,P\,$ corresponds
to the tangent bundle of $\,S^2$.  Since the tangent bundle is known to
be indecomposable, $P$ is also indecomposable (and in particular
$\,P\ncong R^2$), so the class of $P$ in $\,{\mathcal P}_2(R)\,$ is not 
in the image of the map $\,i_1\,$ in (2.3).  This shows that the condition
$\,r\geq d\,$ for the {\it surjective\/} part in (2.6) is the best possible.  
On the other hand, the fact that $\,R\oplus P\cong R\oplus R^2\,$ and 
$\,P\ncong R^2\,$ shows that the condition $\,r\geq d+1\,$ for the
{\it injective\/} part in (2.6) is also the best possible.

\bs
In the noncommutative setting, a very simple case for application is 
that of a group ring $\,R={\mathbb Z}\,\pi$, where $\,\pi\,$ is a finite group.
Here we take $A$ to be $\,{\mathbb Z}\,$ in (2.4), and set $\,d=1$.  It 
follows from (2.4)(1) that any f.g.~projective $\,R$-module
is a direct sum of rank 1 projective modules.  But, according to a famous
example of Swan [Sw$_2$], if $\,\pi\,$ is the generalized quaternion group 
of order 32, there exists a nonfree rank 1 f.g.~projective $\,R$-module $P$ 
such that $\,R\oplus P\cong R \oplus R$.  Since $\,P\ncong R$, we see again 
that, in the conclusion in (2.4)(2), the condition \,rank$\,(P)>d\,$ cannot be 
further weakened.  In Swan's construction, the module $P$ was, in fact, 
chosen such that, for a suitable maximal order $\,S\,$ containing 
$\,R={\mathbb Z}\,\pi$, $\,S\otimes_R P\ncong S$. Therefore, by tensoring up
to $\,S$, we obtain a similar example (of non-cancellation) over a left 
and right hereditary module-finite algebra over the one-dimensional 
ring $\,{\mathbb Z}$.  

\bs
Of course, the bounds $\,r\geq d\,$ and $\,r\geq d+1\,$ are just 
{\it general\/} bounds for surjectivity and injectivity to hold, respectively,
in (2.6).  For specific rings, there can be stronger results.  For
instance, in contrast with Swan's example mentioned above, if $\,\pi_0\,$
is the (ordinary) quaternion group of order 8, then full cancellation holds 
for f.g.~projective modules over the group ring $\,{\mathbb Z}\,\pi_0$, 
according to a result of J.~Martinet [Ma].

\bs
Bass's fundamental results (2.4) and (2.6) in the early 60s provided the 
general framework for much of the subsequent investigations on the
Splitting and Cancellation Problems for projective (and more general) modules 
over various classes of rings.  Without giving any details, let us just 
mention, in this direction, the work of Chase, Mohan Kumar, Murthy, Nori,
Sridharan, Suslin, Swan, Towber, Wiegand, and others on the splitting 
and cancellation
of f.g.~projective modules over affine algebras.  A survey of some aspects of 
this appears in Murthy's article [Mu] in this volume.  On the noncommutative 
side, there are the important cancellation results of Jacobinski, Guralnick,
Levy, Ro\u{i}ter, Swan and many others for lattices over group rings and 
orders in separable algebras.  

\bs
Now it may be said that a cancellation result such as (2.4)(2) is not quite
truly in the noncommutative spirit. In this theorem, the ring $R$ in 
question is a module-finite algebra over a commutative ring $A$; orders in 
finite-dimensional separable algebras are also of the same nature.  Such 
rings are simply {\it not\/} sufficiently representative of a general 
noncommutative ring.  Later developments show that there are indeed 
some ``truly noncommutative'' cancellation theorems, where the cancellation
of a module $\,M\,$ in (2.5) depends on the ``stable range'' of the 
endomorphism ring of $M$ and on the structure of the module $P$.  Such
results, which in essence generalize (2.4)(2), were first found by Warfield 
[Wa$_2$] in 1980 (and in part by Evans in 1973). We shall return to formulate 
these more general cancellation results after we introduce the notion of 
stable range of (noncommutative) rings in \S8.

\bs
\sk
\begin{center}
{\bf \S3.  \ Work Related to Serre's Conjecture} 
\end{center}

\sk
In the Serre-Swan analogy between vector bundles and projective modules, the
counterpart of the affine $\,n$-space $\,k^n\,$ over a field $\,k\,$
is the polynomial ring in $n$ variables $\,R=k[x_1,\dots,x_n]$, so vector
bundles over $\,k^n\,$ correspond to f.g.~projective modules
over $R$.  From the viewpoint of topology, the (real) affine $n$-space is
contractible, so the vector bundles over it are all trivial.  This led 
Serre to ask the question, in his famous FAC paper [S$_1$], {\it whether 
every f.g.~projective module over the polynomial ring $\,R=k[x_1,\dots,x_n]\,$
is free\/} (for any field $k$).  An affirmative answer to this question 
seemed so plausible and convincing to the mathematical public that, almost 
from the very beginning, it became known under the misnomer of ``Serre's 
Conjecture''. This ``conjecture'' is clearly true when $\,n=1$, since in this 
case $\,R=k[x_1]\,$ is a PID, and f.g.~projective modules over a PID are well 
known to be free.  In the general case, it is known (essentially from 
Hilbert's Syzygy Theorem) that any f.g.~projective module $P$ over 
$\,R=k[x_1,\dots,x_n]\,$ is stably free, so Serre's Conjecture boils down 
to a cancellation statement: that $\,P\oplus R\,$ free should imply $P$ 
free.\footnote{In this form, the Conjecture is capable of a completely 
elementary statement meaningful to any student with a high school background 
in algebra: {\it whenever $\,f_1g_1+\cdots + f_rg_r=1\in R$, there is an 
$r\times r$ matrix of determinant $1$ over $R$ with first row 
$\,(f_1,\dots, f_r)$.} See [La$_2$:~(I.4.5)] for details.}  By Bass's 
Theorem (2.6)(2), this would follow, for instance, if the rank of $P$ 
is at least $\,n+1$.

\bs
In the 1960s, Serre's Conjecture became one of the premier open problems in
algebra.  The fact that the Conjecture was prompted by a natural analogy 
with vector bundle theory gave it a certain sense of inevitability; on the
other hand, the fact that the conjecture can be stated so directly and 
in such completely elementary terms made it enticing to all.  Many
algebraists in the 1960s, junior and senior alike, must have tried their 
hands at solving this famous conjecture.  Bass's interest in the structure 
of projective modules, evident from his first two papers [B$_1$, B$_2$], 
naturally steered him in this direction.

\bs
In 1958, Seshadri [Se$_1$] confirmed Serre's Conjecture in the case of two 
variables, proving, as he put it, the ``triviality of vector bundles over 
the affine space $K^2$''.  As a matter of fact, Seshadri showed more
generally that any f.g.~projective module $P$ over $\,A[t]\,$ ``comes from'' 
$A$ if $A$ is a PID [Se$_1$], or the coordinate ring of a nonsingular affine 
curve over an algebraically closed field [Se$_2$].  These results led Bass to 
consider the same problem over $\,R=A[t]\,$ when $A$ is a Dedekind domain. 
The desired goal in this case would be the same as Seshadri's --- to prove 
that any f.g.~projective $P$ over $R$ ``comes from'' $A$, that is, 
$\,P\cong R \otimes_A P_0\,$ for some (necessarily f.g.~projective) module 
$P_0$ over $A$.

\bs
In the early 60s, Bass succeeded in extending Seshadri's argument, and proved
the following result in [B$_4$:~(2.4)].

\bs\nt
{\bf (3.1) Theorem.} {\it Let $A$ be a Dedekind ring, $R=A[t]$, and $P$ be 
a f.g.~$R$-module such that, for any prime ideal $\,{\mathfrak p}\subset A$,
$\,P/{\mathfrak p}P\,$ is a torsionfree $\,(A/{\mathfrak p})[t]$-module.  
Then $P$ is extended from a f.g.~projective $A$-module; in particular, 
$P$ is projective.}

\bs
This implies, in particular, that any f.g.~projective $\,R$-module is extended
from $\,A$.  As it turned out, the same result was obtained independently
by Serre at about the same time; see [S$_3$].  In retrospect, this result 
is perhaps most appropriately called the Seshadri-Bass-Serre Theorem.  

\bs
In 1964, Bass took up what may be considered a noncommutative version 
of Serre's Problem.  Perhaps not too surprisingly, the noncommutative 
case turned out to be more tractable.  Generalizing  the work of 
P.~M.~Cohn, Bass succeeded in proving that, if $\,\pi\,$ is a free group 
(resp.~a free monoid), then for any principal ideal domain $A$, any 
f.g.~projective module over the group ring (resp.~monoid ring) $\,A\,\pi\,$
is free.  This result appeared in Bass's paper [B$_7$], 
in the first volume of the Journal of Algebra.  In the case when $\,\pi\,$ 
is a free monoid on one generator, it retrieves, of course, Seshadri's 
theorem in [Se$_1$].

\bs
Although Bass did not publish further results on Serre's Conjecture after
the 1960s, his keen interest in it continued well into the 70s. When Suslin 
and Vaserstein began to make significant progress on the conjecture in the 
early 70s, the former Soviet Union was still quite isolated mathematically 
from Europe and from the U.S.  In order to make their latest findings known 
to the West, Suslin and Vaserstein could only communicate them by letter to 
Bass.  I still remember vividly the AMS Annual Meeting in San Francisco in 
1974, in which Bass gave an ``impromptu'' lecture on the most recent 
Suslin-Vaserstein results on Serre's Conjecture --- to a room-full of people 
eager to find out how close Serre's Conjecture had come to being solved.  
One of these ``From Russia, with Love'' results proclaimed the freeness of 
projective $\,k[x_1,\dots, x_n]$-modules ``of rank $\geq 1+n/2$''.  This 
surely looked wonderful, but a high-school algebra question unwittingly 
came up: when Suslin wrote ``$1+n/2\,$'' in his letter to Bass, did he mean 
$\,1+\frac{n}{2}$, or could he have meant the more ambitious $\,(1+n)/2\,$?  
It was anybody's guess $\dots\,$.  (As it turned out, Suslin did mean 
$\,1+\frac{n}{2}\,$, as he, perhaps, should.)  Later in June that year, 
in Paris, Bass was to give a similar lecture on the status of Serre's 
Conjecture in the ``S\'eminaire Bourbaki''.  The write-up of these survey 
lectures subsequently appeared in Bass's article [B$_{12}$], with the 
charming title ``Lib\'eration des modules projectifs $\dots \,$''.

\bs
Serre's Conjecture stood open for over twenty years, and was finally proved
in 1976, completely independently and almost simultaneously, by D.~Quillen
[Qu] and A.~Suslin [Su].  (For a detailed exposition on this, see [La$_2$].)  
As is often the case in mathematics, however, the solution of one important 
conjecture was only to be followed by the formulation of a new, more powerful,
conjecture.  After 1976, Serre's Conjecture was generalized into the 
so-called {\it Bass-Quillen Conjecture,} which states the following:

\bs
\nt {\bf (3.2)$_d$} {\it If $A$ is a commutative regular ring\/}\footnote{A 
commutative noetherian ring $\,A\,$ is said to be regular if 
$\,A_{\mathfrak p}\,$ is a regular local ring for any prime ideal 
$\,{\mathfrak p}\subset A$.} {\it of Krull dimension $\leq d<\infty\,$, 
then any f.g.~projective $\,A[x_1,\dots,x_n]$-module is extended from $\,A$.}

\bs
When $A$ is a field, of course, this gives back the original Serre Conjecture
(now the Quillen-Suslin Theorem).  When $A$ is a Dedekind ring and $n=1$, 
$(3.2)_1$ is the Theorem of Seshadri, Bass, and Serre.  By a powerful general 
technique known as ``Quillen Induction'', which is gleaned from Quillen's 
solution of Serre's Conjecture (see [La$_2$:~p.\,139]), one can reduce the 
proof of $(3.2)_d$ to a demonstration of the following special case of it:

\bs\nt
{\bf (3.3)$_d$} {\it If $A$ is a regular local ring of Krull dimension 
$\leq d$, then any f.g.~projective $\,A[x_1]$-module is free.}

\bs
The best result on $(3.2)_d$ and $(3.3)_d$ known to me is that they are 
both true when $\,d\leq 2$, or when $A$ is a formal power series ring over a
field; see [La$_2$:~p.\,138]. On several occasions, I have heard research 
announcements claiming the general truth of $(3.2)_d$ (for all $n$ and all 
$d$), but so far I have not seen any published proofs.  Thus, it appears
that ``Serre's Conjecture'' first raised in the 1950s is still very much 
alive today: it has simply undergone a mathematical metamorphosis and has 
now become the even more challenging ``Bass-Quillen Conjecture''.  Given 
this, we can say with a reasonable amount of certainly that the work of 
Bass on Serre's Conjecture and its generalizations will continue to have 
its impact on the mathematics of the next century.

\bs
\sk
\begin{center}
{\bf \S4.  \ Rings with Binary Generated Ideals: Bass Rings}
\end{center}

\sk
One way to try to prove Serre's Conjecture over $\,R=k[x_1,\dots,x_n]\,$
($k\,$ a field) would be to show that any f.g.~projective $R$-module is
isomorphic to a direct sum of ideals in $R$ (and then use the fact that
Pic$(R)=\{1\}\,$ for the unique factorization domain $R$). The conclusion of 
Bass's Theorem (3.1) has a rather similar flavor: for any Dedekind ring 
$\,A$, this result says that any f.g.~module $\,P\,$ over $\,R=A[t]\,$ 
satisfying the torsionfree hypothesis in that theorem is extended from a 
f.g.~projective $\,A$-module $\,P_0$.  Over the Dedekind ring $\,A$, $\,P_0\,$
is isomorphic to a direct sum of ideals, so $\,P\,$ is likewise isomorphic 
to a direct sum of ideals in $\,R$.  Considerations such as this led Bass 
to the following general question on the decomposition of torsionfree modules:

\bs\nt
{\bf (4.1)} {\it When is it true that any f.g.~torsionfree module 
over a (commutative) noetherian domain $\,R\,$ is isomorphic to a direct 
sum of ideals?} 

\bs\nt
Equivalently, when is it true that any f.g.~indecomposable torsionfree 
$\,R$-module has rank 1 (that is, isomorphic to an ideal of $R$)\,?

\bs
As it turned out, this interesting question led Bass to a fruitful program 
of research.  In his paper [B$_4$], Bass not only proved the result (3.1), 
but also obtained a criterion for the decomposability of all f.g.~torsionfree 
modules into rank one modules over a noetherian domain $R$, under a mild 
assumption on the integral closure $\,\tilde{R}\,$ of $R$.  The main 
theorem (1.7) in [B$_4$] gives the following definitive result.

\bs\nt
{\bf (4.2) Theorem.} {\it For any commutative noetherian domain $R$ such 
that $\,\tilde{R}\,$ is f.g.~as an $\,R$-module, the following two
conditions are equivalent:}

\sk\nt
(1) {\it Any f.g.~torsionfree $\,R$-module is isomorphic to a direct sum 
of ideals\,;}

\nt (2) {\it Any ideal in $R$ can be generated by two elements.}

\bs
In Bass's proof of this theorem, the hypothesis on $\,\tilde{R}\,$ is
needed only for the implication $(2)\Rightarrow (1)$.\footnote{In other
words, $(1)\Rightarrow (2)$ holds for any commutative noetherian domain $R$.}
This hypothesis is, of course, a very natural one from the viewpoint of 
algebraic geometry.  As a matter of fact, Bass's proof (in [B$_4$]) for 
$\,(2)\Rightarrow (1)\,$ in the above theorem is based on a rather 
subtle induction on the length of the $\,R$-module $\,\tilde{R}/R$.  
The fact that this module has finite length is a consequence of the first 
conclusion in the following result describing some of the key properties 
of commutative domains with binary generated ideals.

\bs\nt
{\bf (4.3) Theorem.} {\it If a commutative domain $R$ has the
property that any ideal in $R$ can be generated by two elements, then:}

\sk\nt
(1) {\it $R$ has Krull dimension $\leq 1\,$;}  \\
(2) {\it every ideal in $R$ is reflexive\/}\footnote{An $R\,$-module $P$ is
said to be {\it reflexive\/} if the natural map from $P$ to its double dual 
$\,P^{**}\,$ is an isomorphism.} {\it as an $R$-module\,}.

\bs
Here, the first conclusion, (1), goes back to I.~S.~Cohen. In fact, Cohen 
has proved already in 1949 that the conclusion $\,K$-dim\,$R\leq 1\,$ will
follow if every ideal of the domain $\,R\,$ can be generated by $\,k\,$
elements for a fixed integer $\,k$  [C:~p.\,37, Cor.~1].  Once we have
$\,K$-dim\,$R\leq 1$, the conclusion (2) follows from Lemma (1.6) of [B$_4$].

\bs 
Of course, the driving force behind all of these conditions in (4.2) and
(4.3) is the basic example of a Dedekind ring $R$.  For such a ring, the 
properties (1), (2) in (4.2) and (4.3) are well-known to every student of 
abstract algebra.  The case of Dedekind rings $R$ is precisely the 
{\it integrally closed\/} case of these results, that is, when 
$\,\tilde{R}=R$.  This is, in fact, the beginning case for Bass's proof 
for $(2)\Longrightarrow (1)$ in (4.2) by induction on 
$\,\mbox{length}_R(\tilde{R}/R)$.  How about the case of non-Dedekind 
rings?  It behooves us to recall the two basic examples mentioned by Bass
in [B$_4$:~p.\,324]:

\begin{itemize}
\item  Let $R$ be the ring $\,\{a+2b\,i:\;a,b\in {\mathbb Z}\}$, as a subring 
of the Dedekind ring of Gaussian integers $\,S$.  Here, $\,\tilde{R}=S$, 
and $\,\mbox{length}_R(\tilde{R}/R)=1$.  There are only two isomorphism 
types of nonzero ideals in $R$, namely, $R$ and $\,2R+2i\,R$, so indeed 
every $\,R$-ideal is binary generated. 

\sk
\item  For another typical example, let $\,I\,$ be the additive submonoid 
of the non-negative integers generated by $2$ and an odd integer $\,n\geq 3$, 
and let $\,R\,$ be the subring of $\,S={\mathbb R}[[t]]\,$ consisting of power 
series $\,\sum_{i\geq 0} \,a_it^i\,$ with $\,a_i=0\,$ for $\,i\notin I$. 
Then $\,\tilde{R}=S$, $\,\mbox{length}_R(\tilde{R}/R)=n-2$, and again 
it is easy to check that every ideal of $\,R\,$ is binary generated.
\end{itemize}

\sk
So far we have quoted two main results ((4.2) and (4.3)) from [B$_4$],
which arose from the consideration of torsionfree modules over commutative
noetherian domains, and have the common theme of binary generated ideals 
in such domains.  For Bass, there was another important motivation
for these results\,!  In fact, (4.3) is very much a part of Bass's research 
program studying (not necessarily commutative) noetherian rings of finite 
injective dimensions over themselves.  This is a very important program 
he started in [B$_3$], and continued in [B$_6$].  Let us now explain the 
connections.

\bs
Classically, among the noetherian rings, those that are (say, left) 
self-injective are the so-called {\it quasi-Frobenius} (QF) {\it rings}. 
These rings of self-injective dimension zero have been extensively studied 
in the ring theory literature.  In generalization of this, Bass sought 
characterizations of noetherian rings $\,R\,$ with 
$\,\mbox{id}({}_RR)<\infty$, where ``id'' stands for the {\it injective 
dimension\/} (of a module).  The case of $\,\mbox{id}({}_RR)\leq 1\,$ was 
successfully characterized by Jans [Ja] and Bass [B$_3$:~(3.3)], as follows.

\bs
\nt {\bf (4.4) Theorem.} {\it For any (left and right) noetherian ring $R$, 
the following are equivalent:}

\sk\nt
(1) $\mbox{id}({}_RR)\leq 1${\it ;}  \\
(2) {\it every f.g.~torsionless right $R$-module\/}\footnote{Yes, there 
is a change of side in this statement\,!  Torsionless modules will be
defined in the paragraph following the statement of the Theorem.} 
{\it is reflexive\,;} \\
(3) {\it If ${}_RP$ is f.g.~projective and ${}_RN$ is f.g.~torsionless,
then any short exact sequence 
$\,0\rightarrow P \rightarrow M \rightarrow N \rightarrow 0\,$ splits.}

\bs
Here, since $R$ is not a domain any more, we do not have a natural notion
of torsionfree modules.  In their place, Bass introduced the notion of a
{\it torsionless\/} module: an $\,R$-module $P$ is said to be torsionless
if the natural map $\,\theta_P\,$ from $P$ to its double-dual $\,P^{**}\,$ 
is a monomorphism (that is, for any nonzero $p\in P$, $f(p)\neq 0$ for some 
$\,f\in P^*$).  If $R$ happens to be a commutative domain, a torsionless 
$R$-module $P$ is easily seen to be torsionfree; the converse does not 
hold in general, but does hold if $\,P\,$ is f.g.

\bs
(We should note in passing that, in stating the above theorem (4.4) in
[B$_3$], Bass had another ``equivalent'' condition: ``Every right ideal
in $R$ is reflexive.''  However, while this condition is implied by those
in (4.4), Bass's argument for the converse contained a gap, as was later
acknowledged in [B$_6$:~p.12].  In (6.2) of [B$_6$], this extra condition is
restored in the commutative case, in the form ``$R$ is Cohen-Macaulay,
and every ideal in $R$ is reflexive.'')

\bs
Let us mention two well-known classes of rings that satisfy the conditions
in (4.4).  First, consider any noetherian left hereditary ring $\,R$.  
For $\,E=E({}_RR)\,$ (the injective hull of $\,{}_RR$), we have an exact 
sequence
$$ 0\longrightarrow R\longrightarrow E\longrightarrow E/R\longrightarrow 0.$$
Here, the quotient module $\,E/R\,$ must be {\it injective,} 
by [La$_5$:~(3.22)]. Therefore, we have $\,\mbox{id}({}_RR)\leq 1$.  Another 
interesting (and important) class of noetherian rings $\,R\,$ satisfying
$\,\mbox{id}({}_RR)\leq 1\,$ is given by the group rings $\,kG$, where 
$\,k\,$ is any Dedekind ring, and $G$ is any finite group.  These rings 
are ``one-step away'' from the quasi-Frobenius rings, in that, if $\,a\,$ 
is any nonzero element in $\,k$, the quotient $\,kG/(a)\cong (k/(a))[G]\,$ 
is a well-known example of a QF ring (see [La$_5$:~Exer.~(15.14)]).

\bs
For a commutative noetherian {\it local\/} ring $R$, the condition
$\,\mbox{id}({}_RR)<\infty\,$ turns out to be one of several equivalent 
conditions defining a (local) {\it Gorenstein ring\/}, and for such a 
ring $R$, we have in fact $\,\mbox{id}(R)=K$-dim$\,R$.  This homological 
characterization of a local Gorenstein ring is close in spirit to the
usual homological characterization of regular local rings: recall that, by 
the theorem of Auslander-Buchsbaum and Serre, regular local rings are 
exactly those noetherian local rings $R$ for which we have 
$\,\mbox{gl.dim}(R)<\infty$, or equivalently, $\,\mbox{gl.dim}(R)=K$-dim$\,R$.

\bs
In the case of $\,K$-dim $R=0$, the local Gorenstein rings are precisely 
the local QF rings.  In the case of $\,K$-dim $R=1$, the earliest 
manifestation of local Gorenstein rings was in the form of localizations of 
plane curves, and more generally, complete intersection curves, as was 
noted by Ap\'ery, Samuel, Gorenstein himself, and Rosenlicht.  We refer 
the reader to Craig Huneke's article 
in this volume for a thorough survey on the history of Gorenstein rings.  
In [Hu], Huneke traced the Gorenstein ring notion from the work of the 
above-named authors to that of Grothendieck and Serre, who defined local 
Gorenstein rings in the context of duality theory, via the use of dualizing 
sheaves.  After Serre observed the connection to rings with finite 
self-injective dimensions, Bass wrote the famous ``Ubiquity'' paper [B$_6$] 
in 1963 to put the whole theory of Gorenstein rings on a firm footing.  
One of the basic things he did in this splendid paper was to give a 
``global'' definition for Gorenstein rings: he called a general commutative 
noetherian ring $R$ Gorenstein if all localizations of $\,R\,$ at prime 
ideals are local Gorenstein rings.  If $\,K\mbox{-dim}\,R$ happens to be 
finite, this was shown to be equivalent to $\,\mbox{id}(R)<\infty$, and 
again, in this case, $\,\mbox{id}(R)=K$-dim$\,R$.  Various other 
characterizations for Gorenstein rings (e.g.~in terms of primary 
decompositions of ideals, multiplicities, etc.) are given in the 
``Fundamental Theorem'' in \S1 of [B$_6$].

\bs
To bring the notion of Gorenstein rings to bear on the question of
decompositions of modules considered in (4.2), Bass redid and generalized
this result in \S7 of [B$_6$].  To develop a theory suitable for 
applications to, say, integral representation theory, Bass considered
now a commutative reduced noetherian ring $\,R\,$ of Krull dimension $1$,
with a total ring of quotients $\,K$.  It is no longer assumed that
$R$ is a domain, but we retain the reasonable assumption that 
$\,\tilde{R}$, the integral closure of $R$ in $K$, is a f.g.~$R$-module.
(Of course, $\tilde{R}\,$ is just a finite direct product of
Dedekind rings.)  The f.g.~torsionless $R$-modules can be seen to be
exactly those f.g. $R$-modules $M$ for which the natural map
$\,M\longrightarrow K\otimes_R M\,$ is an {\it injection.}  These 
$R$-modules may be called $\,R$-{\it lattices}, in analogy with the 
terminology used in integral representation theory.  With this setting 
in place, Bass considered the following three conditions on $\,R$:

\bs
(1) {\it Every $\,R$-ideal is generated by two elements.}  

\mk
(2) {\it Any ring between $\,R\,$ and $\,\tilde{R}\,$ is a 
Gorenstein ring.}  

\mk
(3) {\it Every indecomposable $R$-lattice is isomorphic to an $R$-ideal.}

\bs\nt
Bass's main result in \S7 of [B$_6$] is that $(1)\Longleftrightarrow
(2) \Longrightarrow (3)$ (and that (3) fails to imply (1) and (2)
``only in a situation that can be analyzed completely''.\footnote{In
analyzing the situations in which (3) fails to imply (1) and (2), 
however, Bass seemed to have overlooked certain cases. A more complete
analysis when $R$ is local was given later by Nazarova and Ro\u{i}ter in [NR]; 
for the general case, see the paper of Haefner and Levy [HL].  The main
difference between (3) and (1) lies in the fact that (1) is a local
property (as Bass had shown in [Ba$_6$:~(7.4)]), while (3) is {\it not\/}
a local property (see [Gr:~\S2]).}  This result constitutes an expansion 
and simplification of (4.2).  The novel feature about this result is 
the emergence of the ``hereditarily Gorenstein'' condition (2).  
In the post-1963 literature, a ring $R$ in the above setting satisfying 
this condition (2) (or equivalently (1)) has been, quite justifiably, 
called a {\it Bass ring.}

\bs
Later, various authors have obtained new characterizations for Bass rings.
For instance, if $R$ and $\,\tilde{R}\,$ are as above, then the Bass 
ring conditions (1), (2) are further shown to be equivalent to each of 
the following:

\bs
(4) (Greither [Gr]) {\it $\tilde{R}\,$ is binary generated as an 
$R$-module.}

\mk
(5) (Levy-Wiegand [LW]) {\it $\tilde{R}/R\,$ is a cyclic $R$-module.}

\mk
(6) (Wiegand [Wi]) {\it Every faithful $\,R$-lattice has a direct
summand isomorphic to a faithful ideal.}

\bs
Yet other characterizations are obtained by Handelman in [Ha$_2$]. Not 
to be outdone by Bass's ``ubiquity'' title, Handelman called his own paper 
``Propinquity of one-dimensional Gorenstein 
rings''.\footnote{Propinquity$\,:=$ the state of being near in space or 
in time. Handelman used this term to refer to the fact that the rings 
between a Bass ring $R$ and its integral closure $\tilde{R}\,$ are rather 
``close'' to one another.}  In the 80s, Bass was one of the communicating 
editors for the {\it Journal of Pure and Applied Algebra\/}; it was perhaps 
not a coincidence, therefore, that the papers of Handelman, Greither, and
Levy-Wiegand all appeared in that same journal\,!

\bs
In classical commutative algebra, there is a well-known Steinitz-Chevalley 
theory for f.g.~modules over Dedekind domains.  In studying Bass rings,
a natural topic to investigate is therefore the classification problem for
$R$-lattices.  This problem has been successfully tackled by Levy and 
Wiegand.  In [LW], they show, in generalization of the Steinitz-Chevalley
theory, that {\it any given $\,R$-lattice $M$ over a Bass ring $R$ is 
determined by its genus},\footnote{Recall that two $R$-modules $M$ and $N$ 
are said to be in the same genus if $\,M_{\mathfrak m}\cong N_{\mathfrak m}\,$
for every maximal ideal $\,{\mathfrak m}\subset R$.} {\it together with the 
class of a faithful ideal $\,c\ell\,(M)\subseteq R\,$ associated with $\,M$.}
Levy and Wiegand have also obtained very interesting cancellation theorems 
for projective lattices over Bass rings, and later, Levy [Le] even 
extended these results from $\,R$-lattices to general f.g.~$R$-modules 
over a specific class of Bass rings.

\bs
Of course, the most classical examples of Bass rings 
are the Dedekind rings (and their finite direct products).  The next 
class of (non-integrally closed) examples are the quadratic orders, 
that is, $\,{\mathbb Z}$-orders in a quadratic number field.  Any nonzero 
ideal in such an order $\,R\,$ is isomorphic to $\,{\mathbb Z}^2$, 
and is therefore binary generated; therefore, $R$ is a Bass order.
For yet another class of examples, consider the integral group ring
$\,R={\mathbb Z}\,G\,$ for a finite abelian group $\,G$.  As we have
indicated before, $\,R\,$ is always Gorenstein, but it need not be Bass.
Bass's result, in this case, shows that (1), (2), (3) above are in fact
{\it equivalent,} and, by the results of Dade and Heller-Reiner, they 
amount to the fact that the order of $\,G\,$ is square-free (see also
[Gr:~Th.~8.1]).  In other words, $|G|$ being square-free is the necessary 
and sufficient condition for the integral group ring $\,{\mathbb Z}\,G\,$ 
to be a Bass ring. 

\bs
To any perceptive reader of the two papers [B$_3$], [B$_6$], there should
be little doubt that one of the objectives Bass had in mind for his theory
of Gorenstein rings and the decompositions of torsionfree modules was
the potential applications to integral representation theory.  Bass's
explicit mention (in [B$_6$: \S7]) of the $\,R={\mathbb Z}\,G\,$ example in 
the paragraph above was a clear indication of his vision in this direction.
As it turned out, the task of carrying out the program of applying 
Gorenstein (and Bass) rings to integral representation theory was to fall 
on the shoulders of researchers on representation modules in the Russian 
School.

\bs
Of course, in integral representation theory, the rings to be considered
are no longer commutative; but in a sense they are sufficiently close 
to commutative rings.  To put ourselves into this new setting,
we start with a Dedekind ring $R$ with quotient field $K$, and consider
an $\,R$-order $\,\Lambda\,$ in a finite-dimensional separable $K$-algebra
$\,A$. In this setting, we shall be exclusively concerned with (say, left)
$\,\Lambda$-{\it lattices\/}; that is, f.g.~left $\,\Lambda$-modules $L$ which 
{\it embed\/} into the $\,A$-module $\,K\otimes_R L\,$ by the natural map).  
Following Curtis and Reiner [CR:~\S37], we say that $\,\Lambda\,$ is a 
{\it Gorenstein order\/} if the left regular module $\,{}_{\Lambda}{\Lambda}\,$
has the following ``weakly injective'' property: for any $\,\Lambda$-lattices 
$\,M\,$ and $\,N$, any $\,\Lambda$-exact sequence 
$$\,0\rightarrow \Lambda \rightarrow M \rightarrow N \rightarrow 0\,$$
splits.\footnote{This is actually just an equivalent way to say that
$\,\mbox{id}({}_{\Lambda}{\Lambda})\leq 1$: one can see this, essentially, 
by applying $(1)\Leftrightarrow (3)$ in (4.4).  It can also be seen that
$\,\Lambda\,$ being a Gorenstein order is a left/right symmetric
property [CR: (37.8)]; hence the omission of any reference to side.}
We then define {\it Bass orders\/} by a hereditary property: $\,\Lambda\,$ 
is said to be a Bass order if every $\,R$-order in $\,A\,$ containing 
$\,\Lambda\,$ is a Gorenstein order.

\bs 
Comparing these definitions with the traditionally well known ones
for hereditary orders and maximal orders, we can easily verify the 
following hierarchy:
$$
\{\mbox{Maximal\;Orders}\} \subset  \{\mbox{Hereditary\;Orders}\} 
\subset \{\mbox{Bass\;Orders}\} \subset \{\mbox{Gorenstein\;Orders}\},$$
where, as indicated, each inclusion is proper. In parallel to the properties 
(1), (2) and (3) in the commutative case, one can now consider the following 
three properties of a given $\,R$-order $\,\Lambda$:

\bs
(a) {\it Every left ideal of $\,\Lambda\,$ is generated by two elements.}

\mk
(b) {\it $\,\Lambda\,$ is a Bass order.}

\mk
(c) {\it Every indecomposable $\,\Lambda$-lattice is isomorphic to a 
left ideal of $\,\Lambda$.}

\bs
Shortly after the appearance of [B$_6$], Russian workers in integral 
representation theory mounted an ambitious program to try to determine 
the exact relationships between the three properties above in the setting
of noncommutative $\,R$-orders. The definitive results were obtained 
around 1966-67.  In [Ro], Ro\u{i}ter showed that $\,\mbox{(a)} \Rightarrow 
\mbox{(b)}$, and in [DKR], Drozd, Kirichenko, and Ro\u{i}ter showed that
$\,\mbox{(b)} \Rightarrow \mbox{(c)}$.  These results are the best
possible, since in general, (b) does not imply (a), nor does (c) imply
(b).  A detailed exposition on the proofs of the implications
$\,\mbox{(a)} \Rightarrow \mbox{(b)} \Rightarrow \mbox{(c)}\,$
can be found in \S37 of the book of Curtis and Reiner [CR], on which
our present discussion is based.

\bs
The results of Ro\u{i}ter, Drozd and Kirichenko are quite deep, involving
rather subtle analysis of the decompositions of $\,\Lambda$-lattices.
But it was clearly the paradigm of the results of Bass in the commutative 
case that had guided the Russians in their work in this phase of integral 
representation theory.  Subsequently, Bass orders were used as a fundamental 
tool in the work of Drozd-Ro\u{i}ter and Drozd-Kirichenko in their approaches 
to the characterization of orders of finite representation type (the
corresponding work on group rings over rings of algebraic integers having
been completed earlier by H.~Jacobinski).

\bs
In another direction, we should mention that the idea of Gorenstein rings 
has also found recent applications in noncommutative algebraic geometry: 
a notion of (noncommutative) ``Auslander-Gorenstein rings'' has been 
introduced and studied by K.~Ajitabh, S.~P.~Smith, and J.~J.~Zhang 
(see [ASZ]).

\bs
Today, in no small measure due to the influence of Bass's paper, Gorenstein 
rings have lived up to their ``ubiquity'' billing, and are widely used in 
number theory, arithmetic and algebraic geometry, commutative (and
noncommutative) algebra, theory of invariants, and combinatorics.  Even the 
work of Andrew Wiles [W] on elliptic curves and modular forms leading to his 
spectacular proof of Fermat's Last Theorem made use of Gorenstein rings at 
several crucial points: in the Appendix, in Ch.~2 (\S1), and then in 
Ch.~3 on the estimates for the Selmer group.  Recall that it was exactly 
Wiles's earlier attempts to use Euler systems for obtaining the upper 
bounds on the Selmer group that had led to the ``fatal flaw'' in his first 
proof of FLT announced in Cambridge on June 23, 1993.  In 1994, finally 
realizing that the Euler system approach was irreparable, Wiles returned 
to his original approach in estimating the Selmer group using ideas from
Iwasawa Theory.  Basically, in [W:~Ch.~3], Wiles needed to show that certain 
minimal Hecke rings, which are known to be Gorenstein rings, are indeed 
complete intersections.  This was eventually accomplished jointly with Richard
Taylor in [TW].  The role played by commutative algebra (and Gorenstein
rings in particular) in Wiles's paper was described vividly in the 
following words of his [W:~p.\,451]:
\begin{quote}
``The turning point in this and indeed in the whole proof came in the 
Spring of 1991.  In searching for a clue from commutative algebra I had
been particularly struck some years earlier by a paper of Kunz.  I had 
already needed to verify that the Hecke rings were Gorenstein in order
to compute the congruences developed in Chapter 2.  $\,\dots\,$  Kunz's
paper suggested the use of an invariant (the $\,\eta$-invariant $\dots\,$) 
which I saw could be used to test for isomorphisms between Gorenstein 
rings.   A different invariant (the $\,{\mathfrak p}/{\mathfrak p}^2$-invariant
$\dots\,$) I had already observed could be used to test for isomorphisms
between complete intersections.  $\,\dots\,$ Not long afterwards I 
realized that, unlikely though it seemed at first, the equality of 
these invariants was actually a criterion for a Gorenstein ring to be 
a complete intersection.''
\end{quote}

For related literature, see also Lenstra's paper [Len] on complete 
intersections and Gorenstein rings.  In this paper, Lenstra sharpened
Wiles's criterion (in the Appendix of [W]) for a finite $\,{\mathcal O}$-free
local Gorenstein algebra $T$ over a complete discrete valuation ring 
$\,{\mathcal O}\,$ to be a complete intersection.  (Lenstra was able to 
remove the Gorenstein assumption on $T$.)

\bs
\begin{center}
{\bf Part II: Ring Theory}   

\bs
\nt {\bf \S5. Semiperfect Rings as Generalizations of Semiprimary Rings}
\end{center}

\sk
In Part II of this paper, we come to Bass's work in noncommutative ring 
theory.  As we have mentioned in the Introduction, Bass's maiden work 
[B$_1$], developed from his Chicago thesis in 1959, is a ring-theoretic 
paper dealing with the homological generalization of semiprimary rings.
This work turned out to be one of the most influential ring theory papers 
written in that period, as can be partly gauged from the following fact.  
In L.~Small's collected reviews [Sm] of ring theory papers published in 
{\it Math.~Reviews\/} in 1940-79, an average paper got at most a few cross 
citations from other reviews, but Bass's paper [B$_1$] managed to pull as 
many as 29!  Usually, a reviewer would only cite a paper in order to 
indicate the source of a crucial topic or an important idea; the fact 
that, in its twenty years of existence, [B$_1$] drew as many 29 cross 
citations from other reviews was almost without parallel in [Sm].

\bs
What makes [B$_1$] a masterpiece was the fact that it wove together many
of the themes in ring theory and homological algebra that were being 
developed at that time.  On the ring theory side, these themes include:
the Krull-Schmidt Theorem (Azumaya version), chain conditions (suitably
restricted), maximal and minimal submodules (existence questions), the
Jacobson radical (nilpotency questions and lifting of idempotents), and 
Nakayama's Lemma (for general, not necessarily f.g., modules).  
On the side of homological algebra (a pretty new subject in 1959), the 
themes include: projective, injective, and flat modules, projective covers, 
``Ext'' and ``Tor'' functors, and all kinds of homological dimensions.  
Quoting from Bass's comments on ``Theorem P'' in [B$_1$], ``this result 
provides one of those gratifying instances in which several ostensibly 
diverse notions were shown to be intimately related.''  Indeed, it seems 
clear in retrospect that it was this remarkable bridge-building role 
played by the various results in [B$_1$] which helped secure it a permanent 
place in the ring theory literature.  Only a few years after the appearance 
of [B$_1$], the key ingredients of the paper (the theory of perfect and 
semiperfect rings) were incorporated into a standard textbook [L], 
Lambek's ``Lectures on Rings and Modules'', ca.~1966. Today, perfect and 
semiperfect rings continue to be used extensively in a wide variety 
of ring-theoretic settings.

\bs
In this and the next two sections, we'll give a report on Bass's paper 
[B$_1$] and its impact on noncommutative ring theory.  In order to keep the 
size of these sections within bounds, however, we shall only survey below 
the first part of [B$_1$] on perfect and semiperfect rings, and will not 
try to cover its second part on the finitistic homological dimensions 
of rings.

\bs
For the reader's convenience, we first recall a couple of basic 
definitions.  In noncommutative ring theory, a ring $R$ is said to be 
{\it semiprimary\/} if its Jacobson radical rad$(R)$ is nilpotent, and the 
quotient ring $\,R/\mbox{rad}(R)\,$ is (artinian) semisimple. Rings of this 
type were well known in classical ring theory, and had been studied in 
part as a viable generalization of one-sided artinian rings by K.~Asano, 
G.~Azumaya and T.~Nakayama, among others.  With the advent of the new 
style of algebra in the 1950s, homological properties of such rings also 
attracted attention, and had been explored, for instance, in some of 
the papers in the ``Nagoya series'' (ca.~1955-56) on the homological 
dimensions of modules and rings.  Bass had the genius to recognize that the 
artinian condition on $\,R/\mbox{rad}(R)\,$ is important in its own right, 
and later (in [B$_9$:~p.\,504]) defined a ring $R$ to be {\it semilocal\/} 
if it satisfies this condition.  In the commutative case, this condition 
amounts to the finiteness of the number of maximal ideals in $R$, so Bass's 
definition of semilocal rings agrees with the usual one in commutative algebra.
For noncommutative rings, however, ``semilocal'' is no longer equivalent to 
the finiteness of the number of (one-sided or two-sided) maximal ideals.  
Ring theorists now know that Bass had chosen the right definition (as well
as the right name!) for an important class of rings.

\bs
Like semiprimary rings, the semiperfect rings and left (right) perfect 
rings introduced in [B$_1$] are special cases of semilocal rings. Bass 
was led to these rings partly by the classical idea of an injective hull 
due to Eckmann and Sch\"opf.  In general, an {\it injective hull\/} of 
a module $M$ is an injective module $\,I\,$ containing $M$ as a 
{\it large\/} (or {\it essential\/}) submodule, in the sense that:
$$
\forall \;\,\mbox{submodule} \;\,S\subseteq I, \;\;\,S\cap M=0
\Longrightarrow S=0.
$$
Eckmann and Sch\"opf showed that (over any ring) an injective hull 
for $M$ always exists, and is unique up to an isomorphism over $M$.  
In [Ei], Eilenberg initiated a notion of a {\it minimal epimorphism\/}
from a projective module to $M$, and used this notion to study minimal
resolutions, homological dimensions and syzygies.  Bass observed that,
by slightly changing Eilenberg's definition, one gets a precise dual of
the notion of an injective hull: according to Bass, a {\it projective 
cover\/} of a module $M$ is a projective module $P$ with an epimorphism 
$f: P\rightarrow M$ such that ker$(f)$ is {\it small\/} 
(or {\it superfluous\/}) in $\,P\,$ in the sense that:
$$
\forall \;\,\mbox{submodule} \,\;S\subseteq P, \;\;\; S+\mbox{ker}(f)=P 
\Longrightarrow S=P.
$$ 
Such a projective cover for $M$ is easily seen to be unique; the only 
problem is that it may not exist.  For instance, in the category of 
$\,{\mathbb Z}$-modules, the only objects with projective covers (in the 
above sense) are the free abelian groups.

\bs
The task of studying the existence of projective covers in module categories 
was taken up by Bass [B$_1$], who defined a ring $R$ to be left perfect 
(resp.~left semiperfect) if every left $R$-module (resp.~cyclic left 
$R$-module) has a projective cover.\footnote{The adjective ``perfect'' is
due to Eilenberg, who used the term ``perfect categories'' in [Ei] for 
categories of modules possessing projective covers.}
The main work to be done was that of characterizing such rings in 
terms of other interesting conditions.

\bs
In this section, we shall focus on the semiperfect case.  Here, Bass's 
major result is the following.

\bs
\pagebreak
\nt {\bf (5.1) Theorem.} {\it For any ring $R$, the following are equivalent:}

\sk\nt
(1) {\it $R$ is left semiperfect;}   

\sk\nt
(2) {\it Every f.g.~left $R$-module has a projective 
cover;\/}\footnote{It is worth noting that, according to later work of
F.~Sandomierski ([Sa$_1$], [Mue]; see also [La$_4$:~(24.3)]), this 
condition (2) is also equivalent to the following: $(2)'$ {\it Every simple 
right $R$-module has a projective cover.} In a same vein, Sandomierski 
proved that a ring $R$ is right perfect iff every every {\it semisimple\/} 
right $R$-module has a projective cover.}   

\sk\nt
(3) {\it $R$ is semilocal, and idempotents in} $\,R/\mbox{rad}(R)\,$ 
{\it can be lifted to idempotents in $R$.} 

\mk\nt
{\it In case $R$ is a commutative ring, the above conditions are also 
equivalent to:}

\mk\nt
(4) {\it $R$ is a finite direct product of commutative local rings.}

\bs
Note that $(1)\Leftrightarrow (3)$ shows, in particular, the somewhat 
surprising left/right symmetry of semiperfect rings.\footnote{In fact, 
the left-right symmetric characterization of semiperfect rings in (5.1)(3) 
is often used as their {\it definition\/} in the literature: see, for instance,
[La$_3$: p.\,346].}  The second condition in (3) was already quite well 
known in ring theory at that time, and had been studied by Kaplansky,
Jacobson, and Zelinsky, among others.  For instance, if rad$(R)$ is a nil 
ideal, this condition is always satisfied (see [L:~p.\,72, Prop.\,1]). Thus, 
semiperfect rings include all semiprimary rings (which in turn include 
all one-sided artinian rings).  It is in this sense that semiperfect rings 
(and the 1-sided perfect rings to be discussed later) are homological 
generalizations of the classically well known semiprimary rings. 

\bs
The equivalence $(1)\Leftrightarrow (4)$ in Theorem 5.1 shows that, 
for commutative rings, the notion of semiperfect rings has essentially 
nothing new to add to the existing theory.  But the interesting case is that 
of {\it noncommutative\/} semiperfect rings.  Since a major source of 
noncommutative rings is the class of endomorphism rings of modules, it is 
significant to ask when such endomorphism rings are semiperfect.  The answer 
to this question is contained in the following result, which is a remarkable 
extension of the well-known classical theorem that the endomorphism ring of 
a module of finite length is always semiprimary.

\bs\nt
{\bf (5.2) Theorem.} {\it Let $M$ be a right module over a ring $S$.
Then the endomorphism ring} $\,\mbox{End}_S(M)\,$ {\it is semiperfect iff
$M$ has a finite Azumaya decomposition, that is, a decomposition
$\,M=M_1\oplus \cdots \oplus M_n\,$ such that each endomorphism ring}
$\,\mbox{End}_S(M_i)\,$ {\it is local.}

\bs
Although this result was not explicitly stated in [Ba$_1$], it can be proved
easily using the techniques of that paper. (An explicit proof can be found 
in [La$_3$:~(23.6)].)  In some sense, (5.2) shows the ``ubiquity'' of 
semiperfect rings.  For instance, for any {\it injective\/} module 
$\,M_S$, if $M$ has finite uniform dimension, then $\,\mbox{End}_S(M)\,$ is 
semiperfect (and conversely).  As a special case of (5.2), we also see that 
a ring $R$ is semiperfect iff the regular module $\,R_R\,$ has a finite 
Azumaya decomposition; or, in terms
of idempotents, {\it iff there is a decomposition  $\,1=e_1+\cdots +e_n\,$
where the $\,e_i\,$'s are mutually orthogonal idempotents with each
$\,e_iRe_i\,$ a local ring.}  The existence of such a decomposition makes
it possible to generalize a considerable amount of the elementary theory
of artinian rings to a semiperfect ring $R$, including, for instance:
the classification of f.g.~projective $R$-modules, the construction 
of projective covers for f.g.~$R$-modules, the definition of a Cartan 
matrix, block decomposition and basic ring for $R$, etc.  These are parts 
of the foundational material for the theory of artinian rings,
first developed (by Brauer, Osima and others) in the context of group 
algebras of finite groups for applications to the theory of modular 
representations.  It is gratifying to see that a large part of this 
well-known classical theory can be carried over {\it verbatim\/} to 
semiperfect rings.

\bs
The ultimate justification for the introduction of the class of semiperfect
rings lies in the fact that, besides semiprimary rings, there are
also many other natural classes of rings which turn out to be semiperfect.
Let us make a list of some such classes below.

\sk
\begin{itemize}
\item If $N_S$ is any module of finite uniform dimension over a ring $S$
and $M$ is its injective hull, then the endomorphism ring \,End$(M_R)\,$ 
is semiperfect.  (This is essentially a rehash of the remark on injective 
modules we made in the last paragraph.)

\sk
\item  If a semilocal ring $R$ is 1-sided self-injective, then 
$\,R\,$ is semiperfect; see [La$_5$:~(13.4)].

\sk
\item Certain semigroup rings called {\it Kupisch rings\/} by H.~Sato
are semiperfect; see [Sat:~(5.1)].

\sk
\item Any right serial ring (that is, a ring $R$ such that $R_R$ is a
direct sum of uniserial modules) is a semiperfect ring; see, e.g.,
[F$_2$:~p.\,81].

\sk
\item It is well-known to researchers in duality theory that if a ring $R$ 
admits a Morita duality into some other ring $S$, then $R$ and $S$ must
be semiperfect rings.  Indeed, if  $U$ is an $\,(S,R)$-bimodule defining 
the duality, then there is a duality between the Serre subcategories
of $\,U$-reflexive right $R$-modules and left $S$-modules.  Since (left)
$S$-modules have injective hulls, it is not difficult to check from the 
above duality that f.g. (right) $R$-modules have projective covers.  Thus, 
$R$, and hence also $S$, must be semiperfect rings. This basic observation 
is due to Barbara Osofsky.  As a special case, it follows that any
cogenerator ring\footnote{A ring $R$ is called a {\it cogenerator ring\/} 
if $R$ is a cogenerator both as a left and as a right $R$-module. A
cogenerator ring is also known as a {\it Morita ring\/} in the literature.}
is a semiperfect ring.

\sk
\item More generally, Sandomierski [Sa$_2$:~p.\,335] has shown that right 
linearly compact rings are semiperfect, and Azumaya [Az$_1$], Osofsky [Os] and 
others have shown that right PF (pseudo-Frobenius) rings\footnote{A ring 
$R$ is called a {\it right\/} PF {\it ring\/} if any faithful right 
$R$-module is a generator.  For a list of many other characterizations of 
such rings, see [La$_5$:~(19.26)]. The cogenerator rings in Footnote (17)
are precisely the 2-sided PF-rings.} are also semiperfect.
\end{itemize}

\sk
The notion of semiperfect rings has also been generalized in several 
directions. The following are two of them.
\begin{itemize}
\item First, prompted by Bass's characterization (5.1)(2) for a 
semiperfect ring, ring-theorists have come up with a slightly more general 
notion:  a ring $R$ is said to be $F$-semiperfect if every {\it finitely 
presented\/} left $R$-module has a projective cover. This terminology is 
due to Oberst and Schneider [OS]: ``$F\,$'' here stands for ``finite''. 
Clearly, every semiperfect ring is $F$-semiperfect.  In parallel to (5.1)(3), 
there is the following characterization of an $F$-semiperfect ring: {\it
$R$ is $F$-semiperfect iff it is semiregular\/}; that is, idempotents of 
$\,R/\mbox{rad}(R)\,$ can be lifted to $R$, and the epimorphic 
image $\,R/\mbox{rad}(R)\,$ of $R$ is a von Neumann regular ring (instead 
of a semisimple ring). (See, e.g. [OS: (1.2)], and [Ni$_3$, Ni$_4$].)
There is a large supply of such rings.  In fact, if $M$ is any 
quasi-injective module,\footnote{A module $M$ is said to be 
{\it quasi-injective\/} if, for any submodule $\,N\subseteq M$, any 
homomorphism from $N$ to $M$ can be extended to an endomorphism of $M$. 
Needless to say, quasi-injective modules include all injective ones.} then 
the endomorphism ring $\,\mbox{End}_R(M)\,$ is always a semiregular
ring.  This fact, first proved by Faith and Utumi [FU], is one of the 
underpinnings for the Findlay-Lambek-Utumi theory of maximal rings of 
quotients; a self-contained proof of it can be found in [La$_5$:~(13.1)].

\sk
\item   In the theory of direct sum decompositions of modules, there is 
an important class of rings called ``exchange rings'' that are formally 
christened by Warfield in [Wa$_1$].  (A ring $R$ is said to be an
exchange ring if $R_R$ satisfies the exchange property introduced in the 
work of Crawley and J\'onsson [CJ]. For more details on (and 
characterizations of) such rings, see [Ni$_4$].)  The following hierarchy 
shows exactly how semiregular rings and exchange rings compare with 
semiperfect rings:
$$ \{\mbox{Semiperfect\;Rings}\}\subset \{\mbox{Semiregular\;Rings}\} 
                        \subset\{\mbox{Exchange\;Rings}\},  $$
where the second inclusion was shown by Warfield ([Wa$_1$:~Th.~3]; see also
[Ni$_4$: Cor.~(2.3)]).  Here, both inclusions are {\it proper}, as indicated.
It turns out that the main difference between semiperfect rings and 
exchange rings lies in a finiteness condition. In fact,  by combining 
Nicholson's results [Ni$_1$: (4.3)] and [Ni$_4$: (1.9)], one sees that
{\it semiperfect rings are just the exchange rings which do not have 
infinite sets of nonzero orthogonal idempotents\/}; this is stated as
(8.4C) in Faith's recent book [Fa$_3$].  Note that this statement is
much easier to prove if we replace the word ``exchange rings'' by
``semiregular rings'', since it is well-known that, whenever idempotents
can be lifted from $\,R/\mbox{rad}(R)\,$ to $R$, any countable set 
of nonzero orthogonal idempotents in $\,R/\mbox{rad}(R)\,$ can be 
lifted to a similar set in $R$ [La$_3$:~(21.25], and that a von Neumann 
regular ring is semisimple if it has no infinite sets of nonzero orthogonal 
idempotents ([Go$_1$:~(2.16)], [La$_5$:\,Exer.\,6.29]).
\end{itemize}

\sk
Exchange rings are worthy of study in ring theory since they generalize
semiperfect rings, and they form a fairly broad class of rings.  Some known 
facts about semiperfect rings turn out to be true for exchange rings;
the proofs of them are sometimes clearer when they are expressed in the
context of exchange rings.  For instance, M\"uller's well-known result 
[Mue] that {\it any projective (right) module over a semiperfect ring $R$ 
is isomorphic to a direct sum $\,\bigoplus_i\,e_iR\,$ (where 
$\,e_i=e_i^2\in R\,$)} generalizes to, and is quite easy to prove over, any 
exchange ring $R$; see [Wa$_1$:~Th.\,1].  Apply this to a local ring $R$ and
you'll retrieve Kaplansky's classical result (mentioned in \S1) that any 
projective $\,R$-module is free.

\bs
Finally, we should mention that the notion of semiperfect rings has also
been successfully extended to a module-theoretic setting: a module $M$ over 
a ring $R$ is said to be semiperfect if every quotient of $M$ has 
a projective cover.  (Thus, $R$ is semiperfect iff the module $\,{}_RR\,$ 
is.)  The theory of semiperfect modules was initiated in the projective case
in Mares [M] and Kasch-Mares [KM], and has been studied further by Nicholson 
[Ni$_2$, Ni$_3$], and Azumaya [Az$_2$, Az$_4$, Az$_5$], among others. 
For a detailed treatment of semiperfect modules in the general case, see the 
textbook of Kasch [Ka].  Ever since [B$_1$] appeared in 1960, semiperfect 
rings (and modules) have been further studied and extensively utilized 
in many papers in ring theory.  To reflect this trend, I have devoted 
the last chapter of my ring theory graduate text [La$_3$] to an introductory 
exposition on the theory of (perfect and) semiperfect rings.

\bs
\begin{center}
\nt {\bf \S6. Perfect Rings and Restricted DCC}
\end{center}

\sk
In this section, we come to Bass's remarkable characterizations of
left perfect rings; recall that these are, by definition, rings all of
whose left modules have projective covers.  Since left perfect rings
are obviously semiperfect, we expect that one of the characterizations
should be a strengthening of the condition (3) in (5.1); this is given
by the condition (2) below.  The other conditions will be commented
upon later.

\bs
\nt {\bf (6.1) Theorem.} {\it For any ring $R$ with Jacobson radical}
$\,J=\mbox{rad}(R)$, {\it the following are equivalent:}

\mk\nt
(1) {\it $\,R$ is left perfect;}   

\sk\nt
(2) {\it $\,R$ is semilocal, and $J$ is left $T$-nilpotent, that is, for
any $\,a_1,a_2,\dots \,\in J$, $\,a_1 a_2\cdots a_n=0\,$ for some $\,n$;}   

\sk\nt
(3) {\it $\,R$ is semilocal, and every nonzero left $\,R$-module has a 
maximal submodule;}  

\sk\nt
(4) {\it \,Every flat left $R$-module is projective;} 

\sk\nt
(5) {\it $\,R$ satisfies\/} DCC {\it on principal right ideals;}

\sk\nt
(6) {\it \,Any right $R$-module satisfies\/} DCC {\it on cyclic submodules;}

\sk\nt
(7) {\it $\,R$ is semilocal, and every nonzero right $R$-module has a 
minimal submodule;}  

\sk\nt
$(7)'$ {\it $R$ has no infinite sets of nonzero orthogonal idempotents, 
and every nonzero right $\,R$-module has a minimal submodule.}

\bs
It is worth pointing out that left perfect rings are interesting mostly
when they are not semiprime.  In fact, as soon as a left perfect ring $R$
is semiprime, then $R$ is semisimple; see, e.g.~[La$_3$:~(10.24)].

\bs
According to Bass, the ``$T\,$'' in ``$T$-nilpotency'' in (2) stands for
``transfinite''.  Note that this condition in (2) comes somewhere between
$J$ being nilpotent and $J$ being nil.\footnote{In particular, the 
characterization (2) implies immediately that semiprimary rings are both 
right and left perfect. In general, however, right perfect and left perfect
are {\it not\/} equivalent conditions, as was already pointed out by Bass.}  
Similar nil-ness conditions have appeared before in the work of Levitzki 
on the Levitzki radical, but Bass was the first one to realize the role of 
$T$-nilpotency in the study of modules and restricted chain conditions 
on 1-sided ideals.  As it turned out, the successful use of the $T$-nilpotency
condition played a pivotal role in the proof of (6.1).  Results such as (1.6)
showed further the efficacy of the $\,T$-nilpotency conditions.

\bs
We shall now comment on the other conditions in Theorem (6.1).  Along
with (1), the condition (4) is of a homological nature; it can be slightly
rephrased as follows.

\bs\nt
$(4)'$ {\it Every left $R$-module has the same flat as projective dimension.}

\bs\nt
Since, in general, flat modules are precisely direct limits of projective 
modules (by the theorem of Lazard and Govorov), (4) can also be stated 
in the following form:

\bs\nt
$(4)''$ {\it Direct limits of projective left $R$-modules are projective.}  

\bs\nt
Bass's proof for $(4)''\Rightarrow (5)$ is based on a very ingenious 
analysis of the left $R$-module with generators $\,x_1,x_2,\dots,\,$
and relations $\,x_n=a_n x_{n+1}\,$ for all $\,n$, where $\,a_1,a_2,\dots \,$
are given elements in $R$.  

\bs
Note that the first four conditions in (6.1) ((1)--(4)) are conditions on 
{\it left\/} $R$-modules, while the last four conditions ((5)--$(7)'$) are 
on {\it right\/} ones\,!  This switch from left modules to right modules, 
albeit not new for Bass (see Footnote (9)), is in fact one of the inherent 
peculiar features of his Theorem (6.1).  Unfortunately, because of this 
unusual switch of sides, Theorem (6.1) is often misquoted in the literature, 
sometimes even in authoritative sources: see, for instance, [KS:~Thm.\,2, 
p.\,57].  In a couple of standard textbooks, Bass's left (resp.~right) 
$T$-nilpotency condition was renamed the ``right (resp.~left) vanishing 
condition''; in another textbook, the author simply switched Bass's 
definitions of left and right $T$-nilpotency --- and did so without even 
alerting the reader to the difference!  All of this evidently further 
compounded the confusion.  We hereby urge all future authors to exercise 
restraint in changing existing definitions, and to check their statements 
very carefully when quoting Bass's Theorem (6.1).

\bs
While clearly (2) and (5) in (6.1) are purely ring-theoretic 
conditions, Bass's proof of $(2)\Longrightarrow (5)$ was routed through 
the homological condition (4) (or its equivalent versions $(4)'$, $(4)''$).  
This raised the challenging question whether it is possible to give a direct 
proof for $(2)\Longrightarrow (5)$ without using any homological 
algebra.\footnote{The reverse implication $(5)\Longrightarrow (2)$ can 
indeed be done directly without invoking homological tools; see, for 
instance, the exposition in [La$_3$:~p.\,369].}  Such a proof was 
eventually found by Rentschler in [Re].

\bs
Let us now discuss in detail the interesting condition (5), which is a 
natural weakening of the usual artinian condition on right ideals.  Here, 
Kaplansky's ideas played a key role.  In his work with Arens on the 
topological representations of algebras in the late 40s, Kaplansky [K$_1$] 
was motivated to consider the condition stipulating the stabilization of 
all chains of the form 
$$
aR \supseteq a^2 R\supseteq a^3 R \supseteq \cdots 
$$
for all elements $\,a\,$ in a ring $R$.  This condition was shown to be
left/right symmetric many years later by Dischinger, and defines the class
of {\it strongly $\,\pi$-regular rings.}  We will not try to justify the
somewhat clumsy terminology here, for this can be done only by making a
digression into some other definitions not essential for our purposes.
Suffice it to say that, in the case of {\it commutative\/} rings $R$, 
strong $\,\pi$-regularity amounts to $R$ having Krull dimension $0$.
And, as examples of noncommutative strongly $\,\pi$-regular rings, 
Kaplansky mentioned the class of algebraic algebras over a field.  

\bs
In an appendix to [K$_1$], Kaplansky pointed out that it would also be 
natural to consider the DCC for {\it all\/} principal right ideals in a ring,
that is, the condition (5) in Theorem 6.1.  Kaplansky noted that the 
Russian mathematician Gertschikoff had considered rings satisfying this 
condition as far back as 1940, and obtained characterizations of such 
rings that are without nonzero nilpotent elements.\footnote{In fact, 
Gertschikoff's characterization worked more generally for rings possibly 
without an identity element.}  In retrospect, this result of Gertschikoff 
was certainly a harbinger for the equivalence of the conditions (2) and 
(5) in Bass's Theorem 6.1. By a rather strange coincidence, 1959-61 
turned out to be the years in which the minimum condition (5) for 
principal right ideals was destined to blossom: during this period, 
along with the publication of [B$_1$] came two papers of Faith [Fa$_1$] 
and three papers of Sz\'asz [Sz] on the same topic.  Bass invented the 
term ``left perfect'' for his rings (after Eilenberg), while Faith used 
the English acronym ``MP-ring'', and Sz\'asz introduced the German 
acronym ``MHR-ring.''  Carl Faith told me that, in his first meeting
with Bass in the early 60s, they compared notes on their respective
works on rings with DCC on 1-sided ideals, but found surprisingly that 
what they did with these rings had almost nothing in common!

\bs
Going beyond Bass's paper [B$_1$], we should mention several significant 
results obtained later by other authors, which were directly inspired by 
Theorem 6.1.  Prior to the publication of [B$_1$], Chase had proved in
his Chicago thesis that, over any semiprimary ring $R$, DCC holds for 
f.g.~submodules of any (left or right) $R$-module.\footnote{For a
self-contained proof of this, see Faith's paper [Fa$_2$:~p.\,189].}
Prompted by this, Bass asked if, in any left perfect ring, f.g.~right 
ideals satisfy the DCC.  This was affirmed in 1969 by Bj\"ork ([Bj$_1$], see 
also [La$_4$:~(23.3)]) who proved the amazing result that, whenever a 
module (over {\it any\/} ring) satisfies DCC on cyclic submodules, it 
also satisfies the DCC on f.g.~submodules. In particular, this extends 
Chase's result to left perfect rings, showing that in the condition (6) 
in Theorem 6.1, the word ``cyclic'' can be replaced by ``f.g.''.  One 
year later, Jonah [Jo] added ``infinitely many'' more equivalent conditions 
to the list in (6.1):

\bs\nt
$(0)_n$ {\it Any left $\,R$-module satisfies\/} ACC {\it on its 
$\,n$-generated submodules.}  

\bs\nt
(Here, $n\,$ is any natural number.)  Remarkably, these new characterization 
of left perfect rings are in terms of {\it ascending\/} chain conditions; 
also, we are now back full circle to {\it left\/} $\,R$-modules! Note that, 
in particular, Jonah's result implies ``Rings with the minimum condition 
for principal right ideals have the maximum condition for principal left 
ideals.''  This was, in fact, the title of Jonah's paper [Jo].

\bs
To continue our discussion on left perfect rings, we take another look 
at the list of equivalent conditions in (6.1).  Since the condition (7)
is equivalent to the ostensibly weaker condition $(7)'$, an obvious
question one can ask is whether (3) can likewise be weakened to:

\bs\nt
$(3)'$ {\it $R$ has no infinite sets of nonzero orthogonal idempotents, 
and every nonzero left $\,R$-module has a maximal submodule.}

\bs\nt
Note that this condition is in a way ``dual'' to the condition $(7)'$
(just as (3) is dual to $(7)$), so it seems tempting to add $(3)'$ to the 
list of equivalent conditions in (6.1).  An open problem proposed in 
[B$_1$] is, indeed, whether $(3)'$ is a characterization for left-perfectness
of a ring $R$.  Some authors have referred to an affirmative answer to this 
as ``Bass's Conjecture'', although in fact Bass merely raised the question
(see [B$_1$:~p.\,471]).  The answer to this question is possibly a bit 
surprising: it is ``yes'' in the commutative case, and ``no'' in the 
general case, as is shown by Koifman [Ko], and partly by Cozzens, Hamsher,
and Renault. To give an idea of how these conclusions were obtained, we 
proceed as follows.

\bs
In view of (3) and $(3)'$, it is of interest to isolate the condition:
``any nonzero left $R$-module has a maximal submodule''.  Let us say 
that $R$ is a {\it left-max ring\/} if this condition is satisfied.  
It is not hard to see that any left-max ring has a left $T$-nilpotent 
Jacobson radical [La$_4$:~(24.6)]. In the commutative case, one has 
the following characterization of a (left)-max ring, due to Hamsher, 
Renault and Koifman; see [La$_4$:~(24.9)].  

\bs
\nt {\bf (6.2) Theorem.} {\it A commutative ring $R$ is (left)-max iff\/}
$\mbox{rad}(R)$ {\it is $\,T$-nilpotent and} $\,R/\mbox{rad}(R)\,$ 
{\it is a von Neumann regular ring.}

\bs
It follows, in particular, from this result that a commutative max-ring
$R$ is semiregular.  If, in addition, $R$ has no infinite sets of
nonzero orthogonal idempotents, then $\,R/\mbox{rad}(R)\,$ must be
semisimple by an argument we gave in \S5.  Since $\,\mbox{rad}(R)\,$ is 
$\,T$-nilpotent, $R$ is perfect, which answers Bass's question 
affirmatively for commutative rings.

\bs
To treat Bass's question in the general case, we use the notion of a
left $V$-ring --- a ring whose simple left modules are injective.
A well-known characterization for such a ring $R$ is that every left
$R$-module has a zero radical ([La$_5$:~(3.75)]); in particular, $R$ 
must be a left-max ring. Using differential algebra, Cozzens [Co] 
constructed remarkable examples of left $V$-domains that are not 
division rings; such domains clearly satisfy $(3)'$ above, but are 
{\it not\/} left perfect rings, thus answering Bass's question in the
negative.  Moreover, Cozzens's rings are simple, principal left/right 
ideal domains.  The emergence of rings of this type served as an
important turning point for the theory of simple noetherian rings; 
see, e.g.~the book of Cozzens and Faith [CF].

\bs
Many other characterizations of left perfect rings are known.  Without
any attempt at completeness, let us mention a few below.  

\sk
\begin{itemize}
\item P.~A.~Griffith, B.~Zimmermann-Huisgen and others have characterized
left perfect rings in terms of conditions similar to (6.1)(4).  In the 
homological theory of modules, an $\,R$-module $P$ is said to be {\it locally 
projective\/} if every $\,R$-epimorphism $\,f:Q\longrightarrow P$ is locally 
split (in the sense that, for every $\,p\in P$, there exists 
$\,g\in \mbox{Hom}_R(P,Q)\,$ such that $\,fg(p)=p$). It is easy to see 
that such a module $P$ must be flat, so the locally projective modules 
form a class between the class of projective modules and that of flat modules.
In [Zi$_2$], Zimmermann-Huisgen showed that {\it a ring $R$ is left perfect 
iff every locally projective left $R$-module is projective.}  This is, 
therefore, a variant of the criterion (4) for left perfect rings.  

\sk
\item In Theorem 10 {\it loc.~cit.}, Zimmermann-Huisgen has also characterized
left perfect rings in terms of $\,\aleph_1$-separable modules; these are 
modules $P$ with the property that any countable subset of $P$ is contained 
in a countably generated direct summand of $P$.  

\sk
\item In the work of Harada-Ishii [HI], Yamagata [Ya], Zimmermann-Huisgen
and Zimmermann [ZZ$_1$], some characterizations of left perfect rings 
are given in terms of the exchange property of projective left $\,R$-modules.

\sk
\item It seems to be a folklore result that a ring $R$ is left perfect
(resp.~semiperfect) iff every left $R$-module $\,P$ (resp.~f.g.~left $R$-module
$\,P$) is ``supplemented'', that is, every submodule of $\,P\,$ has an
addition complement. (See, e.g.~[Wis:~(42.6), (43.9)].)  Recently,
Keskin [Ke] has further extended these characterizations by using the
more general notion of $\oplus$-supplemented modules.

\sk
\item Optimally, characterizations for left perfect rings should be 
generalizable to characterizations for the endomorphism ring of a given module
to be left perfect.  Some results in this direction can be found in 
[AF$_2$:~\S29], [Wis:~(43.10)], and [Az$_3$].  
\end{itemize}

\sk
In closing, we should mention the fact that much of Bass's work on left
perfect rings can be extended to {\it functor categories.}  In general, if 
{\bf A} is a small additive category, we may view {\bf A} as a generalization 
of a ring (it is a ``ring with several objects''), and even more 
significantly, we may view Add({\bf A}, {\bf Ab}) (the category
of additive functors from {\bf A} into the category of abelian groups)
as a generalization of the category of modules over a ring.
In view of this, it is not surprising that various module-theoretic notions
can be generalized to notions concerning functors in Add({\bf A}, {\bf Ab}).
As it turns out, after defining flat functors, projective covers of functors, 
and DCC for f.g.~subfunctors, etc., one can completely ``transfer'' Bass's
Theorem (6.1) into a theorem on the functor category Add({\bf A}, {\bf Ab}).
This leads then to the notion of a {\it left (resp.~right) perfect\/} 
additive category {\bf A}.  For a full account of all this (as well 
as various other relations between perfect rings and model-theoretic 
algebra), see the book of Jensen and Lenzing [JL].

\bs
Now all of this is not just generalization for generalization's sake!  For 
instance, one can apply it back to the case when {\bf A} is a certain 
additive subcategory of a module category over some ring $\,R$.   By doing 
so, one is sometimes able to make interesting connections, and even
prove nontrivial results.  For instance, if {\bf A} is the category of all
finitely presented modules over a ring $R$, it turns out that the study of 
the perfectness of {\bf A} leads to various insights about (and 
characterizations of) the so-called pure semisimple rings.  For a detailed 
formulation of this, see [JL:~Thm.~B.14].  As a matter of fact, we shall
return to the theme of pure semisimple rings in the next section in the 
context of Auslander's representation theory of artinian rings.

\bs
Needless to say, none of the work discussed above would have been possible 
without the pioneering effort of Bass in [B$_1$].

\bs

\sk
\begin{center}
{\bf \S7. Perfect Rings and Representation Theory}
\end{center}

\mk
In mathematics, a good notion or a good theorem often has a way of finding
surprising connections to things to which it might have seemed unrelated
at first.  As it turned out, Bass's notion of perfect rings and his 
various results on them provide such an example.  The unexpected connections
are to the representation theory of artinian rings, and the study of the
general decomposition theory of modules into direct sums.  In this section, 
we'll give a short account on some of these connections; our discussion
will culminate in an open question in noncommutative ring theory which has 
remained unanswered to this date.  I thank N.~V.~Dung for suggesting
that I include such a discussion here, and for explaining to me the main
results and references in this area of research.

\bs
To set the stage for our discussion, we first recall the following basic 
notion in representation theory.  A ring $R$ is said to be of {\it finite 
representation type\/} (FRT) if it is left artinian and there are only 
finitely many (isomorphism types of) f.g.~indecomposable left $R$-modules.  
According to a result of Eisenbud and Griffith [EG], this notion is 
left/right symmetric, so we are justified in suppressing the word ``left'' 
or ``right'' in referring to rings of FRT.  The study of these rings, in 
the form of finite-dimensional algebras over fields, goes back a long way 
to Brauer, Thrall, Kasch, Kneser, Kupisch, and others.  It turns out that 
rings of FRT have a rather subtle relationship to Bass's perfect rings, 
which we'll now try to explain.

\bs
We start with a notion introduced by Anderson and Fuller [AF$_1$] in the 
decomposition theory of modules. A decomposition of an $\,R$-module $M$ 
into $\,\bigoplus_{i\in I} M_i\,$ ($M_i\neq 0$) is said to ``complement 
direct summands'' if, for each direct summand $\,N\,$ of $M$, there exists 
a subset $\,J\subseteq I$ such that $\,M=N\oplus \bigoplus_{j\in J}M_j$. 
(Note that if such a decomposition for $M$ exists, the summands $\,M_i\,$ 
are necessarily indecomposable.)  Using this terminology and going beyond
[B$_1$], Anderson and Fuller obtained in [AF$_1$] (ca.~1972) the following 
new characterization of left perfect (and semiperfect) rings:

\bs\nt
{\bf (7.1) Theorem.} {\it A ring $R$ is left perfect (resp.~semiperfect) iff 
every projective left $R$-module (resp. f.g.~projective left $R$-module) has 
a direct sum decomposition that complements direct summands, iff the free
left $\,R$-module $\,R\oplus R \oplus \cdots \,$ (resp. $R\oplus R\,$) has
a direct sum decomposition that complements direct summands.}  

\bs
In view of this theorem, it is natural to ask when does {\it every\/} left
$\,R$-module admit a decomposition that complements direct summands.  It
turns out that the answer to this question also involves left perfect rings, 
albeit in different way.  Given a ring $\,R$, let ${}_RU\,$ be the direct
sum $\,\bigoplus_i U_i$, where the $\,U_i\,$'s consist of one isomorphic copy
of each finitely presented left $R$-module.  Now let $E$ be the subring 
(without identity)
of $\,\mbox{End}_R(U)\,$ consisting of endomorphisms $\,f\,$ such that
$\,U_if=0\,$ for almost all $\,i$; $E$ is called the ``left functor ring'' 
of $R$ (for reasons that we shall not elaborate on here).  In 1976, building 
upon the results of Auslander and Harada, Fuller [Ful] proved the following 
remarkable result (see also Simson's paper [Si$_1$]).

\bs\nt
{\bf (7.2) Theorem.} {\it For any ring $R$, the following are equivalent:}

\sk\nt
(1) {\it Every left $\,R$-module has a decomposition that complements direct
summands;} \\
(2) {\it Every left $\,R$-module is a direct sum of f.g.~submodules;}\\
(3) {\it The left functor ring $E$ associated to $R$ (defined above) is left 
perfect.}

\sk\nt
{\it According to a theorem of Chase, these conditions (specifically\/}
(2){\it ) imply that $R$ is a left artinian ring.}

\bs
In (3) above, of course, we'll need to use the notion of a left perfect ring 
without an identity. This is not a big problem; in fact, with essentially
the same definition of left perfectness, Harada [H] has shown that much of 
Theorem 6.1 can be proved for rings without $1$ but ``with enough idempotents''
(see also [Wis:~\S49]).  The condition (2) in (7.2) expresses a very desirable 
module-theoretic property (arbitrary left modules can be ``constructed'' 
from f.g.~ones) that has also been studied by Auslander [Au$_2$], Gruson and 
Jensen [GJ] (and many others in the commutative case).  In [GJ], it was shown 
that (2) is equivalent to $R$ having ``pure left global dimension zero'', 
that is, if every pure short exact sequence of left $\,R$-modules 
splits.\footnote{A short exact sequence of left $\,R$-modules is said to be 
{\it pure\/} if it remains exact upon tensoring by any right $\,R$-module.}   
In 1977, Simson [Si$_2$] introduced the shorter term ``left pure semisimple'' 
for rings $R$ with this property, and later, expanding on the work of Chase 
and Warfield, Zimmermann-Huisgen [Zi$_2$] showed that the property (2) is 
also equivalent to:

\bs\nt
(4) {\it Every left $\,R$-module is a direct sum of indecomposable
submodules.}

\bs\nt
Therefore, {\it any\/} of (1), (2), (3), (4) is a characterization for the 
left pure semisimplicity of a ring $\,R$.

\bs
Now finally, we come to the connection between pure semisimplicity and 
rings of FRT.  This is given by the following theorem of Fuller and Reiten
[FR] (ca.~1975), which completed the earlier work of Auslander, Ringel 
and Tachikawa:

\bs\nt
{\bf (7.3) Theorem.} {\it A ring $R$ is of\/} FRT {\it iff it is both left 
and right pure semisimple; that is, iff every left and every right 
$\,R$-module is a direct sum of \,f.g.~submodules.}

\bs
The major question that remains in this area of study is whether {\it left\/}
pure semisimplicity is equivalent to {\it right\/} pure semisimplicity.  
An affirmative answer to this question is known as the ``Pure Semisimplicity
Conjecture'' (PSC) in ring theory.  If this Conjecture holds up, then the 
sufficiency part of the theorem above would say that left pure semisimple 
rings are of FRT.  As a positive evidence for this, we mention the 
``sparsity'' result of Prest [Pr], Zimmermann-Huisgen and Zimmermann [ZZ$_2$],
which states that, over a left pure semisimple ring $\,R$, there exist only 
finitely many (isomorphism types of) indecomposable left $\,R$-modules 
{\it of any given composition length.}  This means $\,R\,$ comes indeed
``reasonably close'' to being of finite representation type.  Also, it is 
worth noting that, according to results of Simson [Si$_3$] and Herzog [Her],  
the general form of PSC would follow as soon as one could prove that a left 
pure semisimple ring is necessarily {\it right\/} artinian.  For more 
perspectives and recent results on PSC, see [Si$_4$] and the literature 
referenced therein.

\bs
For connections with the work of Bass, we note that, with the help of left 
functor rings, PSC can actually be formulated entirely in terms of the 
notion of perfect rings.  We need the following additional fact, which can 
be gleaned from the work of Auslander [Au$_1$] and Fuller [Ful] (see also 
Wisbauer's book [Wis:~(54.3)]).

\bs\nt
{\bf (7.4) Proposition.} {\it A ring $R$ is left and right pure semisimple
iff its left functor ring $\,E\,$ (defined in the paragraph preceding\/}
(7.2){\it ) is both left and right perfect.}

\bs
Granted this and (7.2), the issue of whether a left pure semisimple ring $R$ 
needs to be right pure semisimple boils down to testing whether the 
perfectness of $\,E\,$ on the left would imply its perfectness on the right. 
In an earlier footnote (Footnote (20)), we have mentioned Bass's observation
that, in general, a left perfect ring need not be right perfect. The above 
considerations show, however, that PSC reduces to a ``left perfect implies 
right perfect'' statement, {\it for the special class of left functor rings.}

\bs
As we have mentioned at the end of \S6, Bass's ideas of using projective 
covers and the condition ``flat $\Rightarrow$ projective'' for module 
categories have been successfully carried over to categories of 
functors.  Such a generalization has proved to be quite fruitful for the 
representation theory of artinian rings; see, {\it inter alia}, the work 
of Auslander [Au$_1$] and Simson [Si$_1$].  Bass's left and right 
$\,T$-nilpotency conditions for ideals perhaps also inspired Auslander's 
definition in [Au$_1$] of noetherian and conoetherian conditions for 
families of homomorphisms between modules.  In fact, for a left artinian 
ring $R$, the left and right $\,T$-nilpotency conditions on the Jacobson 
radical of the left functor ring $\,E\,$ are precisely equivalent to the 
noetherian and conoetherian conditions for families of homomorphisms between 
f.g.~indecomposable left $\,R$-modules.  From these and the preceding 
discussions, we see that Bass's early work in [B$_1$] has certainly played 
a substantial and rather interesting role in the later development of the 
representation theory of artinian rings.

\bs

\sk
\begin{center}
{\bf \S8. Stable Range of Rings}
\end{center}

\mk
The notion of stable range for rings is another great invention of Bass
that has proved to be of lasting importance in algebra and ring theory.  This 
notion was originally introduced by Bass (in [B$_9$:~p.\,14]) for the study 
of stabilization questions in algebraic $K$-theory, or at least for working 
with the functor $\,K_1$.  The definition of stable range goes as follows.

\bs\nt
{\bf (8.1) Definition.}  
We say that a positive integer $\,n\,$ is {\it in the stable range\/} 
of a ring $R$ (or, more informally, that $R$ ``has stable range 
$\,n\,$'') if, whenever $\,a_1R+\cdots +a_{m+1}R=R$ for $m \geq n$ (where all 
$a_i\in R$), there exist elements $\,x_1,\dots, x_{m} \in R\,$ such that 
$$\,(a_1+a_{m+1}x_1)R+\cdots +(a_{m}+a_{m+1}x_{m})R=R.  $$

\sk
It is straightforward to see that, to verify the condition for stable 
range $n$ given above, it is sufficient to do so in the critical case 
$\,m=n$.  (It follows that if $n$ is in the stable range for $R$, then 
so is any larger integer.)  Also, although it appears that we should have 
referred to the above condition as $R$ having ``right'' stable range $n$ 
(since the definition was based on the use of right ideals), it has been 
shown later by Vaserstein [V$_1$: Th.~2] (and also by Warfield 
[Wa$_2$:~Th.~(1.6)]) that ``right stable range $n$'' and 
``left stable range $n$'' are actually {\it equivalent\/} conditions.  
For this reason, we shall suppress any reference to side in referring to 
the stable range conditions defined in (8.1).

\bs
For readers familiar with Bass's big book [B$_{10}$], we should point out
that there is a small discrepancy between the stable range notation used
here and that used in [B$_{10}$].  What we called ``stable range $n$'' above
corresponds to what Bass called SR$_{n+1}$ in his book.  We believe
the usage in (8.1) is now the standard one.

\bs
Bass's study of the stable range of rings was again motivated by the basic 
ideas of stability in the homotopy theory of vector bundles.  His main results
may be summarized in the fundamental theorem (8.2) below, which is to be
thought of as the ``$K_1$-analogue'' of Corollary (2.6).
Here, $\mbox{GL}_n(R)\,$ denotes the $n\times n$ general linear group over the 
ring $R$, and $\,\mbox{E}_n(R)\,$ denotes its subgroup generated by the 
$\,n\times n\,$ elementary matrices $\,I_n+ae_{ij}$ ($i\neq j,\;a\in R$); 
$K_1(R)\,$ denotes the Whitehead group  $\,\mbox{GL}(R)/E(R)\,$ in algebraic 
$K$-theory, where $\mbox{GL}(R)$ and E$(R)$ are, respectively, the direct
limits of the groups $\mbox{GL}_r(R)\,$'s and $\mbox{E}_r(R)$, taken over 
all positive integers $\,r$.

\bs\nt
{\bf (8.2) Theorem.} (See [B$_{9}$:\;(4.2), (11.1)] and 
[B$_{10}$:\;pp.\,239-240].)  

\sk\nt
(1) {\it If a ring $R$ has stable range $n$, then}
$\,\mbox{GL}_n(R)\rightarrow K_1(R)\,$ {\it is a surjective homomorphism.  
Moreover, for any $\,r\geq n+1$}, $\mbox{E}_r(R)\,$ {\it is normal in} 
$\,\mbox{GL}_r(R)$, {\it and} 
$\mbox{GL}_r(R)/\mbox{E}_r(R)\rightarrow K_1(R)\,$ {\it is an isomorphism.} 

\sk\nt
(2)  {\it Let $A$ be a commutative ring whose maximal 
ideal spectrum} $\,\mbox{max}(A)\,$ {\it is a noetherian space of dimension 
$\, d$.  Then, any module-finite 
$\,A$-algebra $R$ has stable range $\,d+1$.  In particular, the conclusions 
of\/} (1) {\it apply to such an algebra $R$ for $\,n=d+1$.}

\bs
The second part of this theorem has been generalized to noncommutative
noetherian rings $R$ by Stafford; see [St]. In the case of commutative
rings $R$, there has also been work done toward the removal of the 
noetherian hypothesis on the maximal ideal spectrum.  For instance, 
Heitmann [Hei] has shown that, if $R$ is a commutative ring of Krull 
dimension $d$, then $R$ has stable range $d+2$ (and $d+1$ if $R$ is a
domain).

\bs
As one would perhaps expect, the case of stable range $1$ has a special
significance, not shared by the general case of stable range $n$. It is true, 
for instance, that having stable range $1$ is a Morita invariant property 
of a ring, while, for $\,n>1$, having stable range $n$ is {\it not\/} a Morita 
invariant property [V$_1$, V$_2$, Wa$_2$].  The study of stable range 
has led to many new and unexpected results in the arithmetic 
of rings and the structure theory of modules.  In the rest of this 
section (and the next), we shall give a survey of some of the interesting 
mathematics which resulted from this study, and which certainly would not 
have been possible without the pioneering work of Bass.  We begin with 
the following useful observation that apparently first appeared in an 
unpublished note of Kaplansky [K$_3$].

\bs\nt
{\bf (8.3) Proposition.}  {\it If a ring $R$ has stable range $1$, then $R$ 
is Dedekind-finite; that is, $\,uv=1\in R$ implies that $\,vu=1\in R$.}

\bs\nt
{\bf Proof} (following [K$_3$]).  Say $\,uv=1\in R$. Since $\,vR+(1-vu)R=R$,
there exists $\,t\in R\,$ such that $\,w:=v+(1-vu)t\,$ has a right inverse.
Left-multiplying this equation by $u$, we get $\,uw=1$. Therefore, $w$ has
also a left inverse.  It follows that $\,w,\,u\in \mbox{U}(R)$ (the group of
units of $R$), and hence $\,vu=1$. \ \ \ \ \ \ \ \ {\bf QED}

\bs
For an alternative approach to this Proposition, see the beginning of
\S8  below.  An immediate consequence of the Proposition is the following 
somewhat sharper formulation for the condition of stable range $1$.

\bs\nt
{\bf (8.4) Corollary.} {\it A ring $R$ has stable range $\,1\,$ iff, whenever 
$\,aR+bR=R$, there exists $\,x\in R\,$ such that} $\,a+bx\in \mbox{U}(R)\,$
{\it (the group of units of $R$).}

\bs
We recorded this corollary explicitly since the characterization for stable 
range 1 contained herein is often used as its {\it definition\/} in papers 
in the literature dealing with the stable range of rings.

\bs
In [B$_9$: (6.5)], Bass proved the following basic result:

\bs\nt
{\bf (8.5) Theorem.} {\it Any semilocal ring $R$ has stable range $1$.}

\bs
For a modern reader, this is a very natural result.  One can first check
(without too much difficulty) that a semisimple ring has stable range $1$, 
and then deduce (8.5) from the observation that, in general, $R$ has $n$ 
in its stable range if $\,R/\mbox{rad}(R)\,$ does. (For more details, see
[La$_3$:~pp.\,313-314].)   Conceptually, one can think of (8.5) as the
``$0$-dimensional case'' of part (2) of Theorem 8.2; in fact, in the case
of commutative rings, $R$ semilocal means that $\,\mbox{max}(R)\,$ is
{\it finite,} and therefore a 0-dimensional space.

\bs
Bass's study of stable range was quickly picked up by Estes and Ohm, who
obtained various results in 1967 on stable range in the commutative case.  
The case of local rings shows already that a (commutative) ring $R$ may have 
arbitrary Krull dimension, but still have $1$ in its stable range.  In [EO], 
Estes and Ohm constructed commutative B\'ezout domains of arbitrary dimension
with $1$ in their stable range.  For other interesting examples of rings
of stable range 1 and 2, see, for instance, the work of Heinzer [He],
Jensen [Je], Menal-Moncasi [MM], and Warfield [Wa$_2$].

\bs
Warfield's paper [Wa$_2$] is a very important contribution to the study of 
the stable range of {\it noncommutative\/} rings.  In this paper, he studied 
the endomorphism ring $\,E:=\mbox{End}_R(M)\,$ of a right module $M$ over 
an arbitrary ring $R$, and sought characterizations for such a ring $E$ 
to have stable range $n$ (for a given integer $n$) . Since Warfield's
results lead to very interesting noncommutative cancellation theorems
harkening back to Bass's results in \S2, we shall give a quick exposition
on the gist of [Wa$_2$] below.  The first result of Warfield 
characterizes the stable range of $\,E:=\mbox{End}_R(M)\,$ in terms of
a certain ``substitution property'' of $M$, as follows.

\bs\nt
{\bf (8.6) Theorem.} [Wa$_2$:~(1.6)] {\it Let $M$ be any right $R$-module
with endomorphism ring} $\,E:=\mbox{End}_R(M)$ {\it (operating on the left
on $M$). Then $\,E\,$ has (right) stable range $n$ iff $\,M\,$  has the 
following ``$\,n$-Substitution Property'':}

\mk\nt
$(S)_n$ {\it For any split epimorphism $\,\pi\,$ from $T=M^n\oplus X\,$
to $M$, there exists a splitting $\,\varphi\,$ such that 
$T=\varphi(M)\oplus Y\oplus X$, where $\,Y\subseteq M^n$.}

\bs
The proof of (8.6) consists of a fairly straightforward manipulation of the 
definition of stable range $\,n$.

\bs\nt
{\bf (8.7) Remarks.} (1) Note that if $(S)_n$ holds, then in the notation 
there we have 
$$  \mbox{ker}(\pi) \cong Y\oplus X, \;\;\,\mbox{and} \;\;\,
M^n\cong Y\oplus M.  \leqno (8.8)  $$
The first follows since both sides of the equation are direct complements
of $\,\varphi(M)\,$ in $T$; the second follows similarly since
both $M^n$ and $\,Y\oplus \varphi(M)\,$ are direct complements of $X$ 
in $T$ (and $\,\varphi(M)\cong M$).  

\sk\nt
(2) Since the $\,n$-Substitution Property $(S)_n$ on $M$ 
depends only on $E$ by (8.6), it follows that if $N$ is any module over 
any other ring with endomorphism ring isomorphic to $E$, then $M$ satisfies 
$(S)_n$ iff $N$ does.

\bs
To see the effect of (8.6) (and (8.7)) on module cancellations, we formulate 
the following slight improvement on Warfield's result in [Wa$_2$:~(1.3)].

\bs\nt
{\bf (8.9) Theorem.}  {\it Let $M$ be any right $R$-module with endomorphism
ring} $\,E:=\mbox{End}_R(M)$ {\it having stable range $n\geq 1$.  For any right
$R$-modules $P$ and $Q$, the following are equivalent:}

\mk\nt
(1) {\it For some module $X$, we have $\,P=M^{n-1}\oplus X$ and 
$\,M\oplus P \cong M\oplus Q\,$;}  \\
(2) {\it For some modules $X$ and $Y$, we have}
$$P\cong M^{n-1}\oplus X, \;\;\, Q\cong Y\oplus X, \;\;\,
\mbox{and} \;\;\, M\oplus Y\cong M^n. \leqno (8.10) $$

\bs\nt
{\bf Proof.}  If (2) holds, then from (8.10):
$$M\oplus P\cong M\oplus M^{n-1}\oplus X\cong M^n\oplus X
\cong M\oplus Y\oplus X \cong M\oplus Q.$$
Conversely, if (1) holds, let $\,T:=M\oplus P=M\oplus M^{n-1}\oplus X$.  
Since $T\cong M\oplus Q$ has a split epimorphism $\,\pi\,$ into $M$ with 
$\,\mbox{ker}(\pi)\cong Q$, the property $(S)_n$ for $M$ implies the 
existence of a splitting $\,\varphi$ such that 
$\,T=\varphi(M)\oplus Y\oplus X$, where $\,Y\subseteq M^n$.  By  (8.8), 
we have $M\oplus Y\cong M^n$, and 
$\,Q\cong \mbox{ker}(\pi) \cong Y\oplus X$, as required in (2).  
(Note that the condition $\,Y\subseteq M^n\,$ from $(S)_n$ is never used
in all of the arguments above.) \ \ \ \ \ \ \ {\bf QED}

\bs
Note that, if we could have the implication
$$\,M\oplus Y\cong M^n \Longrightarrow Y\cong M^{n-1}, \leqno (8.11) $$
then in (2) above, we would have been able to conclude that $\,P\cong Q$.  
While (8.11) does hold sometimes,\footnote{By a principle of A.~Dress 
[La$_5$:~(18.59)], (8.11) holds iff $\,E\oplus W\cong E^n\,$ (as right 
$E$-modules) implies $\,W\cong E^{n-1}$.  Thus, (8.11) will hold if, say, 
stably free right modules over $E$ are free and $E$ has invariant basis 
number.} it cannot be counted on in the most general situation.  To get 
a good cancellation theorem out of (8.9), we can, instead, impose a 
slightly stronger assumption on the module $\,P$, as in the following 
result.

\bs\nt
{\bf (8.12) Warfield's Cancellation Theorem.} [Wa$_2$:~(1.2)]
{\it Let $M$ be any right $R$-module with endomorphism ring} 
$\,E:=\mbox{End}_R(M)$ {\it having stable range $n\geq 1$. If $\,P\,$ 
contains a direct summand $\,M^n$, then for any $\,R$-module $Q$,}
$$ M\oplus P \cong M\oplus Q \Longrightarrow P\cong Q. $$

\bs\nt
{\bf Proof.} Write $\,P=M^n\oplus W=M^{n-1}\oplus X$, where
$\,X:=M\oplus W$. If $\, M\oplus P \cong M\oplus Q$, we'll have the
ismorphisms in (8.10) for some $Y$, and hence
$$ Q\cong Y\oplus X \cong Y\oplus M \oplus W \cong M^n\oplus W = P,$$
as desired. \ \ \ \ \ \ \ \ {\bf QED}

\bs
In the special case when $\,E=\mbox{End}_R(M)\,$ has stable range 1, 
the theorem above holds even {\it without any assumption\/} on $\,P\,$
(or on $\,Q$); we shall come back to this point in \S9.

\bs
In comparison with Bass's Cancellation Theorem (2.4)(2), the advantage of 
Warfield's (8.12) lies in the fact that there is neither a dimension 
assumption nor a noetherian assumption in its statement. Thus, (8.12) may
be regarded as a ``truly noncommutative'' cancellation result. To see that 
(8.12) essentially retrieves the cancellability of a f.g.~projective module 
$\,M\,$ in the setting of (2.4)(2), we may first replace $M$ there by 
$\,R^m\,$ for some $\,m$, and then reduce to the case when $\,M=R$.  In this 
case, $\,E=\mbox{End}_R(M)\cong R$. If $R$ is module-finite over a 
commutative ring whose maximum spectrum is a noetherian space of dimension 
$\,d$, then by (8.2)(2), $R$ has stable range $\,d+1$. So, as long as 
$\,P\,$ has a free direct summand of rank $d+1$, (8.12) enables us to 
cancel off $\,M=R\,$ in (2.5).

\bs 
In order to apply (8.12), one needs to know about the stable range of the 
endomorphism rings of modules.  In this direction, Warfield has further
generalized Bass's result (2.4)(2) by taking the module-finite $\,A$-algebra 
$R$ there and considering finitely presented right modules $M$ over $R$.  
In [Wa$_2$:~(3.4)], Warfield showed that, if max$(A)$ is noetherian of 
dimension $\,d$, then $\,\mbox{End}_R(M)\,$ has $\,d+1\,$ in its stable 
range.  Thus, one can apply (8.12) with $\,n=d+1$.  Since $\,M_R\,$ need 
no longer be a projective module, this leads to cancellation results well 
beyond the reach of (2.4)(2).  Warfield's proof involved some new techniques,
since Bass's methods did not apply to endomorphism rings.

\bs

\sk
\begin{center}
{\bf \S9. Rings of Stable Range One}
\end{center}

\mk
To conclude our exposition, we shall consider in this section the case
of rings of stable range 1.  The crucial result here is that of Bass 
(Theorem 8.5), which states that semilocal rings have stable range 1.
Again, as it turned out, this single result served as the fountain-head 
of many beautiful ideas to come; the survey in this section will show 
how much this result has stimulated subsequent research on the stable 
range of rings.

\bs
For the balance of this section, we shall only consider the stable 
range 1 case.  First let us make the following useful observation.

\bs\nt
{\bf (9.1) Remark.} {\it Any module $M_R$ satisfying the 1-Substitution 
Property $(S)_1$ in} (8.6) {\it is necessarily Dedekind-finite; that is, 
$M$ is not isomorphic to a proper direct summand of itself.}  

\bs
To see this, suppose there is an isomorphism $\,\pi: M\oplus W \rightarrow M$. 
Applying $(S)_1$ to $\,\pi$, we have a splitting $\,\varphi\,$ for $\,\pi\,$ 
(necessarily $\,\pi^{-1}\,$!) such that 
$\,M\oplus W=\varphi(M)\oplus Y\oplus W\,$ for some $Y$. But
$\,\varphi(M)=T$, so $W$ (as well as $Y$) must be zero, proving (9.1).

\bs
Note that (9.1) gives us another view of the fact (8.3) that a ring 
$R$ of stable range 1 is Dedekind-finite.  In fact, if $R$ has (right) 
stable range 1, then, viewing $R$ as $\,\mbox{End}(R_R)$, (8.6) shows that
$\,R_R\,$ satisfies $(S)_1$, and hence (9.1) shows that $R_R$ is
Dedekind-finite.  But this is the same as saying that its endomorphism
ring $R$ is Dedekind-finite (as a ring).  This is a more conceptual, if 
somewhat longer, proof of (8.3).

\bs
Equipped now with the information in (9.1), let us take another look 
at the condition $\,(S)_1\,$ on $\,M$:

\bs\nt
$(S)_1$ {\it For any split epimorphism $\,\pi\,$ from $T=M\oplus X\,$
to $M$, there exists a splitting $\,\varphi\,$ such that 
$T=\varphi(M)\oplus Y\oplus X$, where $\,Y\subseteq M$.}

\bs\nt
Here, we have $M\cong Y\oplus M$ (as observed in (8.8)), so if $(S)_1$
holds for $\,M$, (9.1) implies that $Y=0$, and hence the conclusion of 
$(S)_1$ simplifies to $\,T=\varphi(M)\oplus X$. Thus, $(S)_1$ simply says 
that, if $X$ and $\,\mbox{ker}(\pi)\,$ both have direct complements isomorphic 
to $M$ in a module $T$, then they have a {\it common\/} direct complement 
in $T$.  We can therefore restate (8.6) in the case $\,n=1\,$ as follows.

\bs\nt
{\bf (9.2) Theorem.}  {\it The endomorphism ring} $\,E:=\mbox{End}_R(M)\,$ 
{\it of an $R$-module $\,M_R\,$ has stable range $1$ iff $\,M\,$  has the 
following ``Substitution Property'': Whenever a right $R$-module $T$ has 
direct decompositions $\,T=M_i\oplus P_i\,$ for $\,i=1,2$ where 
$\,M_1\cong M_2\cong M$, there exists a submodule $\,C\subseteq T\,$ such 
that $\,T=C\oplus P_i\,$ for $\,i=1,2$. } 

\bs
Note that, in the above notation, if the submodule $C$ exists, we have, in 
particular, $\,P_1\cong T/C\cong P_2$.  Therefore, we have the following
consequences of (9.2).

\bs\nt
{\bf (9.3) Corollary.} (1) {\it If an $R$-module $\,M_R\,$ has an endomorphism 
ring with $1$ in its stable range, then $\,M\,$ is ``cancellable'';
that is, for any right $\,R$-modules $\,P,\;Q$,}
$$  M\oplus P\cong M\oplus Q  \,\Longrightarrow \, P\cong Q. $$
\nt (2) {\it If a ring $R$ has stable range $1$, then $R_R$ is cancellable,
and hence so is any f.g.~projective right (respectively, left) 
$R$-module.}

\bs
The result (9.3) was a cancellation theorem obtained earlier (ca.~1973) by 
E.~G.~Evans in [Ev].  We are reporting the results in the reverse chronological
order here only because Warfield's result (9.2) is more general than Evans's, 
so it is logically more convenient to state (9.2) first and deduce (9.3) 
as its corollary.  

\bs
To put things in the right historical perspective, we should 
also point out that what we called the ``Substitution Property'' in (9.2)
had been considered as early as 1971 by L.~Fuchs.  In [Fu], apparently 
unaware of Bass's work on stable range, Fuchs obtained a different 
characterization of the Substitution Property, and proved the result (9.2) 
in the case $\,M_R=R$ (and some other cases), falling just short of proving 
(9.2) for {\it any\/} $R$-module $\,M_R$.  Fuchs also considered the case 
of von Neumann regular rings $R$ and sought conditions for $R$ to have the 
Substitution Property (that is, to have stable range $1$). He obtained
the following result ([Fu: Cor.\,1 and Th.\,4]), which was also proved 
independently by Kaplansky [K$_3$: Th.\,3].

\bs\nt
{\bf (9.4) Theorem.}  {\it A von Neumann regular ring 
$R$ has stable range $1$ iff, for $\,x,\,y\in R$, $\,xR\cong yR\Longrightarrow 
R/xR\cong R/yR$.}

\bs
Since any principal right ideal $\,xR\,$ in a von Neumann regular ring $R$
is a direct summand of the right module $R_R$ (and $\,R/xR\,$ is isomorphic 
to any direct complement of $\,xR$), the latter condition in Theorem 9.4
amounts to an ``Internal Cancellation Property'' of the module $R_R$; that 
is, isomorphic direct summands in $R_R$ have isomorphic direct complements.
Now by subsequent results of Ehrlich and Handelman on regular rings $R$ ([Eh],
[Ha$_1$]; see also [La$_4$:~p.\,242]), the module $R_R$ has this property 
iff the ring $R$ is {\it unit-regular}, in the sense that, for any $\,a\in R$,
there exists a unit $\,u\in R\,$ such that $\,a=aua$.  Therefore, (9.4)
amounts to the fact that {\it a von Neumann regular ring $R$ has stable 
range $1$ iff $R$ is unit-regular\/}.  This statement seemed to have appeared
explicitly for the first time in Henriksen's paper [Hen:~Prop.\,8].

\bs
While the above results on unit-regular rings are concerned mainly with
the cancellation of f.g.~{\it projective\/} modules, we should mention at 
least one case where (9.3)(1) applies more generally to {\it all\/} 
f.g.~modules. This is a result due to Menal [Me], which deals with a 
certain subclass of unit-regular rings.

\bs\nt
{\bf (9.5) Theorem.} {\it Let $R$ be a von Neumann regular ring all of
whose primitive factor rings are artinian.  Then for any f.g.~projective
right $R$-module $M$,} $\,\mbox{End}_R(M)\,$ {\it has stable range} $1$.
{\it In particular, for any $R$-modules $P$ and $Q$, 
$\,M\oplus P\cong M\oplus Q\,$ implies $\,P\cong Q$.}

\bs
Now let us come back to Bass's theorem (8.5). Combining this theorem 
with Evans's Cancellation Theorem (9.3), 
we see that {\it any module with a semilocal endomorphism ring is 
cancellable.}  Another recent result of Facchini, Herbera, Levy and V\'amos
[FH:~(2.1)] showed that, if $P$, $Q$ are modules with semilocal endomorphism 
rings, then they have the ``$n$-cancellation property'' (for any
integer $n\geq 1$); that is,  
$$  P^n\cong Q^n \Longrightarrow P\cong Q.  $$  
(For a detailed survey on $\,n$-cancellation, see [La$_6$].)  These 
remarkable results give us strong motivation for finding classes of modules 
which have semilocal endomorphism rings.  Now the problem of describing the 
endomorphism rings of specific kinds of modules has had a time-honored 
history, starting with the very famous Schur's Lemma:

\bs\nt
{\bf (9.6)}  If $M$ is a simple $R$-module, $\,\mbox{End}_R(M)\,$ 
is a division (and hence semilocal) ring.

\mk\nt
{\bf (9.7)} Also well-known is the following classical generalization of 
Schur's Lemma: {\it if $\,M_R\,$ is a module of finite length, 
then} $\,\mbox{End}_R(M)\,$ {\it is a semiprimary (and hence semilocal) 
ring\/}; see, for instance, [La$_4$: Ex.~(21.24)].

\mk\nt
{\bf (9.8)}  Any strongly indecomposable module has, by definition,
a local (and hence semilocal) endomorphism ring.  Examples include, for
instance, all indecomposable modules of finite length, and all injective 
(in fact all quasi-injective or pure-injective) indecomposable modules; 
see, e.g. [La$_5$: (3.52), Exer.\;(6.32)] and [F$_2$:~(2.27)].

\mk\nt
{\bf (9.9)}  Camps and Dicks [CD] showed in 1993 that {\it any artinian
module has a semilocal endomorphism ring.}

\mk\nt
{\bf (9.10)} In 1995, Herbera and Shamsuddin [HS] generalized the
Camps-Dicks result above by showing that, {\it if a module $M$ has finite 
uniform dimension and co-uniform dimension\/},\footnote{A module $M$ is said 
to have finite uniform dimension if there is a bound on the numbers $\,n\,$
for which $M$ contains a direct sum of $n$ nonzero modules; dually, $M$ is 
said to have finite co-uniform dimension if there is a bound on the numbers
$n$ for which $M$ has a quotient that is a direct sum of $n$ nonzero modules.}
{\it then it has a semilocal endomorphism ring.}  This class includes all 
linearly compact modules\footnote{An $\,R$-module $\,M\,$ is said to be 
{\it linearly compact\/} if, for any submodules $\,N_i\subseteq M\,$ and 
any elements $\,m_i\in M$, any system of congruences 
$\,\{x\equiv m_i\; (\mbox{mod}\;N_i)\}\,$ that is finitely solvable is 
solvable.}, and therefore all artinian modules.  Thus, the Herbera-Shamsuddin
result implies that linearly compact modules have semilocal endomorphism
rings; this affirms an earlier conjecture of Faith. In the {\it commutative\/}
case, a stronger conclusion is possible.  In this case, any linearly compact
module $\,M_R\,$ is algebraically compact (i.e.\,purely injective) 
according to Jensen and Lenzing [JL:\,p.\,289].  Therefore, combining the 
Herbera-Shamsuddin result above with [JL:\,(7.5)], we see that $\,M\,$ 
has a {\it semiperfect\/} endomorphism ring.  This provides a nice 
connection back to the material of \S5.

\ms\nt
{\bf (9.11)} A nonzero module is said to be {\it uniserial\/} if its 
submodules form a chain under inclusion.  Such a module has clearly
uniform dimension and co-uniform dimension both equal to $1$. Thus,
if $\,M:=U_1\oplus \cdots \oplus U_n\,$ where the $\,U_i\,$'s are
uniserial, then the ``serial module'' $\,M\,$ has uniform and co-uniform 
dimension $\,n$, and hence it has a semilocal endomorphism ring by (9.10).
For more information on the structure of $\,\mbox{End}(U_i)$,
$\,\mbox{End}(M)$, and its applications to a (weak) Krull-Schmidt Theorem
for serial modules, see Facchini's paper [F$_1$] and his recent book
[F$_2$].

\bs
In view of the remarks made before (9.6), any modules of the type listed in 
(9.6)--(9.11) above have both the cancellation property and the 
$\,n$-cancellation property.

\bs
New results on stable range $1$ are still being discovered today. Let us just
mention a couple of more recent works. In [Ca], Canfell studied the stable 
range $1$ condition from the viewpoint of completing diagrams of modules 
using automorphisms, and made connections to the notions of epi-projective 
and mono-injective modules.  In [Ar], inspired by earlier results of Menal
and Goodearl-Menal [GM], Ara proved the remarkable result that

\bs\nt
{\bf (9.12) Theorem.} {\it Any strongly $\,\pi$-regular ring} 
({\it see} \S6) {\it has stable range $1$.}  

\bs
In the commutative case, this amounts to the fact that any ring of Krull 
dimension $0$ has stable range 1, which is close in spirit to part (2) of 
Bass's Theorem (8.2).  In the general case, Ara's result is rather deep, and 
depends heavily on the use of noncommutative techniques. Combined with
the earlier result of Armendariz, Fisher and Snider [AFS], (9.12)
implies that any module with the ``Fitting decomposition property'' is
cancellable: see (5.4) in [La$_6$].

\bs
In closing, we should also mention some interesting variations of the
stable range one condition, due to Goodearl and Menal.  In [GM], Goodearl and
Menal considered the following two conditions for a ring $\,E$:

\bs\nt
(A) $\,\forall \;a,b\in E, \;\;aE+bE=E \,  \Longrightarrow   \,  
\exists \;\,u\in \mbox{U}(E)\,$ such that $\,a+bu \in \mbox{U}(E)$; and

\mk\nt
(B) $\,\forall \;a,b\in E, \; \exists \;u\in \mbox{U}(E) \,$
such that $\, a-u \;\,\mbox{and} \,\; b-u^{-1} \in \mbox{U}(E).$

\bs\nt
They showed that $(B)\Rightarrow (A)$, and obviously 
$(A) \Rightarrow $ stable range $1$.  They called condition (A) 
(right) ``unit 1-stable range''; the unitary analogue of it has been
exploited for $C^*$-algebras.  Another variation, also due to Goodearl, is a 
{\it weakening\/} of stable range 1 condition. Extending Bass's definition, 
Goodearl defined a ring $\,E\,$ to have the {\it power-substitution 
property\/} if, whenever $\,aE+bE=E$, there exists a positive integer $\,n\,$ 
and a matrix $\,X\in {\mathbb M}_n(E)\,$ such that $\,a\,I_n+b\,X\,$ is 
a unit in $\,{\mathbb M}_n(E)$.   (This property is also known to be 
left-right symmetric.)   Clearly, if $\,E\,$ has stable range 1, then it has 
the power-substitution property.  The converse of this is, however, not true.  
For instance, it can be shown that the ring $\,{\mathbb Z}\,$ has the 
power-substitution property, but it certainly does not have stable range 1.
The {\it raison d'\^etre\/} for the  power-substitution property
lies in the following important result of Goodearl on ``power-cancellation''
[Go$_2$: Cor.~4]:

\bs
\nt {\bf (9.13) Theorem.} {\it Let $\,M,\,P,\,Q\,$ be right modules over
an arbitrary ring $\,R\,$ such that $\,M\oplus P\cong M\oplus Q$.
If\/} $\,E:=\mbox{End}_R(M)\,$ {\it has the power-substitution property,
then $\,P^n\cong Q^n\,$ for some $\,n \geq 1$.}

\bs
This result of Goodearl is, of course, an extension of (9.3)(1).
In fact, in the case when the endomorphism ring $E$ has stable range $1$, 
the proof of Goodearl's result boils down to that of Evans, and gives the 
conclusion $\,P\cong Q$.  The power of (9.13) stems from the fact that many 
different types of rings happen to have the power-substitution property;
a good list of such rings is given in my survey [La$_6$:~p.\,34].  In \S5 
of this survey, there is also a summary of the many interesting results, due 
to Goodearl, Guralnick, and Levy-Wiegand, on the power-cancellation exponent
$\,n\,$ occurring in the statement of (9.13).

\bs

\sk\nt
{\bf Acknowledgment.} The author heartily thanks Larry Levy for a critical
reading of a penultimate version of this paper.   His many thoughtful 
comments and suggestions have led to various improvements in my exposition.

\mk


\begin{thebibliography}{ABC}
\vspace{.3cm}

\bibitem[AF$_1$]{Ander1}  F.~W.~Anderson and K.~R.~Fuller: {\it Modules with
decompositions that complement direct summands}, J.~Algebra {\bf 22}(1972),
241-253.

\bibitem[AF$_2$]{Ander2}  F.~W.~Anderson and K.~R.~Fuller: {\it Rings and 
Categories of Modules},  Second Edition, Graduate Texts in Math., 
Vol.~{\bf 13}, Springer-Verlag, Berlin-Heidelberg-New York, 1992.

\bibitem[ASZ]{ASZ} K.~Ajitabh, S.~P.~Smith and J.~J.~Zhang: 
{\it Auslander-Gorenstein rings}, preprint, 1998.

\bibitem[Ar]{Ar} P.~Ara: {\it Strongly $\,\pi$-regular rings have stable
range $1$}, Proc.~Amer.~Math.~Soc. {\bf 124}(1996), 3293-3298.

\bibitem[AFS]{Afs} E.~P.~Armendariz, J.~W.~Fisher and R.~L.~Snider: {\it On
injective and surjective endomorphisms of finitely generated modules},
Comm. Algebra {\bf 6}(1978), 659-672.

\bibitem[Au$_1$]{Au1} M.~Auslander: {\it Representation theory of artin 
algebras II\/}, Comm.~Algebra {\bf 1}(1974), 269-310.

\bibitem[Au$_2$]{Au2} M.~Auslander: {\it Large modules over artin algebras}, 
in "Algebra, Topology and Categories", Academic Press, 1976, pp.\,1-17. 

\bibitem[Az$_1$]{Az1} G.~Azumaya: {\it Completely faithful modules and 
self-injective rings,} Nagoya Math. J. {\bf 27}(1966), 697-708.

\bibitem[Az$_2$]{Az2} G.~Azumaya: {\it Characterization of semiperfect and 
perfect modules}, Math.~Zeit. {\bf 140}(1974), 95-103.

\bibitem[Az$_3$]{Az3} G.~Azumaya: {\it Locally split submodules and modules
with perfect endomorphism rings}, in ``Noncommutative Ring Theory'' 
(Athens, Ohio), Lecture Notes in Mathematics, Vol.~{\bf 1448}, 
Springer-Verlag, Berlin-Heidelberg-New York, 1990.

\bibitem[Az$_4$]{Az4} G.~Azumaya: {\it $F$-semiperfect modules},
J.~Algebra {\bf 136}(1991), 73-85.

\bibitem[Az$_5$]{Az5} G.~Azumaya: {\it A characterization of semiperfect 
rings and modules}, in {\it Ring Theory,} Proc.~Ohio State-Denison
Conf.~(1992) (S.K.~Jain and S.~T.~Rizvi, eds.), World Scientific
Publishers, River Edge, N.J., 1993.

\bibitem[B$_1$]{Ba1}   H.~Bass: {\it Finitistic dimension and a homological
generalization of semiprimary rings}, Trans. Amer. Math. Soc. {\bf 95}(1960), 
466-488.

\bibitem[B$_2$]{Ba2}   H.~Bass: {\it Projective modules over algebras},
Annals of Math. {\bf 73}(1961), 532-542.

\bibitem[B$_3$]{Ba3} H.~Bass: {\it Injective dimension in Noetherian rings},
Trans. Amer. Math. Soc. {\bf 102}(1962), 18-29.

\bibitem[B$_4$]{Ba4}   H.~Bass: {\it Torsion free and projective modules},
Trans. Amer. Math. Soc. {\bf 102}(1962), 319-327.

\bibitem[B$_5$]{Ba5} H.~Bass: {\it Big projective modules are free}, 
Illinois J.~Math. {\bf 7}(1963), 24-31.

\bibitem[B$_6$]{Ba6} H.~Bass: {\it On the ubiquity of Gorenstein rings},
Math. Zeit.~{\bf 82}(1963), 8-28.

\bibitem[B$_7$]{Ba7}   H.~Bass: {\it Projective modules over free groups 
are free}, J. Algebra {\bf 1}(1964), 367-373.

\bibitem[B$_8$]{Ba8}   H.~Bass: {\it The stable structure of quite general
linear groups}, Bull. Amer. Math. Soc. {\bf 70}(1964), 429-433.

\bibitem[B$_9$]{Ba9}   H.~Bass: {\it $K$-theory and stable algebra}, 
Publ.~IHES {\bf 22}(1964), 5-60.

\bibitem[B$_{10}$]{Ba10}   H.~Bass: {\it Algebraic $\,K$-Theory}, Mathematics
Lecture Note Series, W.~A.~Benjamin, Inc./Addison-Wesley, Reading, MA., 1968.

\bibitem[B$_{11}$]{Ba11}   H.~Bass: {\it Lib\'eration des modules projectifs
sur certains anneaux des polyn\^omes}, S\'em. Bourbaki, 26e ann\'ee (1973/74),
Exp. No. 448, pp.\,228-254, Lecture Notes in Math., Vol.\,{\bf 431}, 
Springer-Verlag, Berlin-Heidelberg-New York.

%\bibitem[B$_{12}$]{Ba12}   H.~Bass: {\it The Morita Theorems}, mimeographed 
%lecture notes, University of Oregon, Eugene, Oregon, 1962.

% bibitem[B$_{13}$]{Ba13}   H.~Bass: {\it Descending chains and the Krull
% ordinal of commutative noetherian rings}, J.~Pure Appl.~Algebra
% {\bf 1}(1971), 347-360.

\bibitem[BMS]{BMS}   H.~Bass, J.~Milnor and J.-P.~Serre: 
{\it Solution of the congruence subgroup problem for $SL_n\;(n\geq 3)$
and $\,Sp_{2n}\;(n\geq 2)$}, Publ.~IHES {\bf 33}(1967), 59-137.

%\bibitem[BG$_1$]{BG1} H.~Bass and R.~Guralnick: {\it Projective modules with
%free multiples and powers}, Proc.~Amer.~Math.~Soc.~{\bf 96}(1986), 207-208.

%\bibitem[BG$_2$]{BG2} H.~Bass and R.~Guralnick: {\it Torsion in the Picard
%group and extension of scalars}, J.~Pure Appl.~Algebra {\bf 52}(1988),
%213-217.

\bibitem[BS]{BS} H.~Bass and S.~Schanuel: {\it The homotopy theory
of projective modules}, Bull. Amer. Math. Soc. {\bf 68}(1962), 425-428.

\bibitem[Bj$_1$]{Bj1}  J.-E.~Bj\"ork: {\it Rings satisfying a minimum condition
on principal ideals}, J. reine angew. Math. {\bf 236}(1969), 112-119.

\bibitem[Bj$_2$]{Bj2}  J.-E.~Bj\"ork: {\it On rings satisfying certain
chain conditions}, J. reine angew. Math. {\bf 237}(1970), 63-73.

\bibitem[CD]{CD} R.~Camps and W.~Dicks: {\it On semilocal rings},
Israel J. Math. {\bf 81}(1993), 203-211.

\bibitem[Ca]{Ca}  M.~J.~Canfell: {\it Completions of diagrams by automorphisms
and Bass' first stable range condition},  J.\;Algebra {\bf 176}(1995). 480-503.

% \bibitem[Ch]{Ch} S.~U.~Chase: {\it Direct products of modules},
% Trans. Amer. Math. Soc. {\bf 97}(1960), 457-473.

\bibitem[C]{C} I.~S.~Cohen: {\it Commutative rings with restricted 
minimum condition}, Duke Math.~J. {\bf 17}(1950), 27-42.

\bibitem[Co]{Co} J.~H.~Cozzens: {\it Homological properties of the ring of
differential polynomials}, Bull. Amer. Math. Soc. {\bf 76}(1970), 75-79.

\bibitem[CF]{CF} J.~H.~Cozzens and C.~Faith: {\it Simple Noetherian Rings},
Cambridge Tracts in Math., Vol.~{\bf 69}, Cambridge Univ.~Press, 
Cambridge, 1975.

\bibitem[CJ]{CJ}  P.~Crawley anf B.~J\'onsson: {\it Refinements for infinite
direct decompositions of algebraic systems}, Pacific J.~Math.~{\bf 14}(1964), 
797-855.

\bibitem[CR]{CR}  C.~W.~Curtis and I.~Reiner: {\it Methods of Representation 
Theory}, Vol.~1, J.~Wiley and Sons, Inc., New York, 1981.

\bibitem[DKR]{DKR}  Yu.~A.~Drozd, V.~V.~Kirichenko and A.~V.~Ro\u{i}ter: 
{\it On hereditary and Bass orders,} Izv.~Akad.~Nauk SSSR, Ser.~Mat.
{\bf 31}(1967), 1415-1436. (English Transl.: Math. USSR Izv. {\bf 1}(1967),
1357-1376.)

\bibitem[Eh]{Eh}  G.~Ehrlich: {\it Units and $1$-sided units in regular 
rings}, Trans. Amer. Math. Soc. {\bf 216}(1976), 81-90.

\bibitem[Ei]{Ei}  S.~Eilenberg: {\it Homological dimension and syzygies},
Annals of Math. {\bf 64}(1956), 328-336.

\bibitem[EG]{EG}  D.~Eisenbud and P.~A.~Griffith: {\it The structure of serial
rings}, Pacific J.~Math. {\bf 36}(1971), 109-121.

\bibitem[EO]{EO}  D.~Estes and J.~Ohm: {\it Stable range in commutative 
rings}, J.~Algebra {\bf 7}(1967), 343-362.

\bibitem[Ev]{Ev}  E.~G.~Evans, Jr.: {\it Krull-Schmidt and cancellation
over local rings}, Pacific J.~Math. {\bf 46}(1973), 115-121.

\bibitem[F$_1$]{F1} A.~Facchini: {\it Krull-Schmidt fails for serial 
modules}, Trans. Amer. Math. Soc. {\bf 348}(1996), 4561-4575.

\bibitem[F$_2$]{F2} A.~Facchini: {\it Module Theory: Endomorphism Rings
and Direct Sum Decompositions in Some Classes of Modules}, Progress in
Math., Vol.\,{\bf 167}, Birk\"auser Boston, 1998.

\bibitem[FH]{FH} A.~Facchini, D.~Herbera, L.~Levy and P.~V\'amos: {\it
Krull-Schmidt fails for artinian modules}, Proc. Amer. Math. 
Soc.~{\bf 123}(1995), 3587-3592.

\bibitem[Fa$_1$]{Fa1}  C.~Faith: {\it  Rings with minimum condition on
principal ideals}, Arch.~Math. {\bf 10}(1959), 327-330; Part II, {\it ibid.}
{\bf 12}(1961), 179-181.

\bibitem[Fa$_2$]{Fa2}  C.~Faith: {\it  Rings with ascending condition on
annihilators}, Nagoya Math.~Journal {\bf 27}(1966), 179-191.

\bibitem[Fa$_3$]{Fa3}  C.~Faith: {\it  Rings and Things and a Fine Array
of Twentieth Century Associative Algebra}, Math.~Surveys and Monographs,
Vol.~65, 1999, Amer. Math. Soc., Providence, R.I.

\bibitem[FU]{FU}  C.~Faith and Y.~Utumi: {\it Quasi-injective modules and
their endomorphism rings}, Arch. Math. {\bf 15}(1964), 166-174.

\bibitem[Fu]{Fu}  L.~Fuchs: {\it On a substitution property of modules}, 
Monatsh. Math. {\bf 75}(1971), 198-204.

\bibitem[Ful]{Ful}  K.~R.~Fuller: {\it On rings whose left modules are direct
sums of finitely generated modules}, Proc. Amer. Math. Soc. {\bf 54}(1976), 
39-44.

\bibitem[FR]{FR}  K.~R.~Fuller and I.~Reiten: {\it Note on rings of finite 
representation type and decompositions of modules}, Proc. Amer. Math. Soc.
{\bf 50}(1975), 92-94.

\bibitem[Go$_1$]{Go1}   K.~R.~Goodearl: {\it von Neumann Regular Rings},
Second Edition, Krieger Publ.~Co., Malabar, Florida, 1991.

\bibitem[Go$_2$]{Go2}  K.~R.~Goodearl: {\it Power-cancellation of modules},
Ring Theory II, Proc. Second Conf. Univ. Oklahoma (B.~R.~McDonald and
R.~A.~Morris, eds.), pp.~131-147, Lecture Notes in Pure and Applied Math. 
{\bf 26}, M.~Dekker, New York, 1977.

\bibitem[GM]{GM}  K.~R.~Goodearl and P.~Menal: {\it Stable range one for
rings with many units}, J.~Pure Appl.~Algebra {\bf 54}(1988), 261-287.

\bibitem[Gr]{Gr}  C.~Greither: {\it On the two generator problem for the
ideals of a one-dimensional ring}, J.~Pure Appl.~Algebra {\bf 24}(1982), 
265-276.

\bibitem[GJ]{GJ}  L.~Gruson and C.~U.~Jensen: {\it Deux applications de la
notion de $L$-dimension}, C. R. Acad. Sci. Paris, S\'er.~A {\bf 282}(1976), 
23-24.

%\bibitem[GJ]{GJ}  L.~Gruson and C.~U.~Jensen: {\it Modules alg\'ebriquement
%compacts et foncteurs\/} lim$^{(i)}\,$, C.~R.~Acad.~Sci. Paris, S\'er.~A
%{\bf 276}(1973), 1651-1653.

\bibitem[HL]{HL}  J.~Haefner and L.~Levy: {\it Commutative orders whose 
lattices are direct sums of ideals}, J.~Pure Appl.~Algebra {\bf 50}(1988), 
1-20.

\bibitem[Ha$_1$]{Ha1}  D.~Handelman: {\it Perspectivity and cancellation in 
regular rings}, J.~Algebra {\bf 48}(1977), 1-16.

\bibitem[Ha$_2$]{Ha2}  D.~Handelman: {\it Propinquity of one-dimensional
Gorenstein rings}, J.~Pure Appl.~Algebra {\bf 24}(1982), 145-150.

\bibitem[H]{H}  M.~Harada: {\it Perfect categories\/} I, II, III,
Osaka J.~Math. {\bf 10}(1973), 329-341, 343-355, 357-367.

\bibitem[HI]{HI}  M.~Harada and T.~Ishii: {\it On perfect rings and the
exchange property}, Osaka J.~Math. {\bf 121}(1975), 483-491.

\bibitem[He]{He}  W.~J.~Heinzer: {\it $J$-Noetherian integral domains with
$1$ in the stable range}, Proc. A.M.S. {\bf 19}(1968), 1369-1372.

\bibitem[Hei]{Hei}  R.~Heitmann: {\it Generating ideals in Pr\"ufer domains},
Pacific J.~Math.~{\bf 62}(1976), 117-126.

\bibitem[Hen]{Hen}  M.~Henriksen: {\it On a class of regular rings that are
elementary divisor rings}, Arch.~Math. {\bf 24}(1973), 133-141.

\bibitem[HS]{HS}  D.~Herbera and A.~Shamsuddin: {\it Modules with semilocal
endomorphism rings}, Proc. Amer. Math. Soc.~{\bf 123}(1995), 3593-3600.

\bibitem[Her]{Her} I.~Herzog: {\it A test for finite representation type},
J.~Pure Appl.~Algebra {\bf 95}(1994), 151-185.

\bibitem[Hu]{Hu}  C.~Huneke: {\it Hyman Bass and ubiquity: Gorenstein rings},
this volume, pp.

% \bibitem[J]{J} N.~Jacobson: {\it Structure of Rings}, A.M.S. Colloquium
% Publications, Vol.\,{\bf 37}, 1956. (Revised Edition: 1964.)

\bibitem[Ja]{Ja} J.~P.~Jans: {\it Duality of Noetherian rings}, Proc. Amer.
Math. Soc. {\bf 12}(1961), 829-835.

\bibitem[Je]{Je} C.~U.~Jensen: {\it Some curiosities of rings of 
analytic functions}, J.~Pure Appl.~Algebra {\bf 38}(1985), 277-283.

\bibitem[JL]{JL} C.~U.~Jensen and H.~Lenzing: {\it Model Theoretic Algebra:
with particular emphasis on fields, rings, and modules}, Gordon and Breach, 
N.Y., 1989.

\bibitem[Jo]{Jo}  D.~Jonah: {\it Rings with the minimum condition for 
principal right ideals have the maximum condition for principal left ideals},
Math. Zeit. {\bf 113}(1970), 106-112.

\bibitem[K$_1$]{K1} I.~Kaplansky: {\it Topological representations of 
algebras\/} II, Trans. Amer. Math. Soc. {\bf 68}(1950),62-75.

\bibitem[K$_2$]{K2} I.~Kaplansky: {\it Projective modules}, Ann. Math.
{\bf 68}(1958), 372-377.

\bibitem[K$_3$]{K3} I.~Kaplansky: {\it Bass's first stable range condition},
mimeographed notes, 1971.

\bibitem[Ka]{Ka} F.~Kasch: {\it Modules and Rings}, Academic Press, 
London-New York, 1982.

\bibitem[KM]{KM} F.~Kasch and E.~Mares: {\it Eine Kennzeichnung 
semi-perfekter Moduln}, Nagoya Math.~J. {\bf 27}(1966), 525-529,

\bibitem[Ke]{Ke} D.~Keskin: {\it Characterization of right perfect rings
by $\,\bigoplus$-supplemented modules}, Lecture at International Conference 
on Algebra and Its Applications, Athens, Ohio, March, 1999.

\bibitem[Ko]{Ko} L.~A.~Koifman: {\it Rings over which each module has a
maximal submodule\/}, Mat.~Zametki {\bf 7}(1970), 359-367.
(English Transl.: Math.~Notes {\bf 7}(1970), 215-219.)

\bibitem[KS]{KS} A.~I.~Kostrikin and I.~R.~Shafarevich: {\it Algebra\/} II,
Encyclopedia of Math.~Sciences, Springer-Verlag, Berlin-Heidelberg-New York,
1991.

\bibitem[La$_1$]{La1}  T.~Y.~Lam: {\it The category of noetherian modules}, 
Proc.~Nat.~Acad.~Sci.~{\bf 55}(1966), 1038-1040.

\bibitem[La$_2$]{La2}  T.~Y.~Lam: {\it Serre's Conjecture}, Lecture
Notes in Math., Vol.\,{\bf 635}, Springer-Verlag, Berlin-Heidelberg-New York,
1978.

\bibitem[La$_3$]{La3}   T.~Y.~Lam: {\it A First Course in Noncommutative 
Rings}, Graduate Texts in Math., Vol.~{\bf 131}, Springer-Verlag, 
Berlin-Heidelberg-New York, 1991.

\bibitem[La$_4$]{La4}  T.~Y.~Lam: {\it Exercises in Classical Ring Theory},
Problem Books in Mathematics, Springer-Verlag, Berlin-Heidelberg-New York,
1995.

\bibitem[La$_5$]{La5} T.~Y.~Lam: {\it Lectures on Modules and Rings},
Graduate Texts in Math.,  Vol.\,{\bf 189}, Springer-Verlag, 
Berlin-Heidelberg-New York, 1999.

\bibitem[La$_6$]{La6} T.~Y.~Lam: {\it Modules with isomorphic multiples
and rings with isomorphic matrix rings --- a survey}, Monographie No.~{\bf 35},
L'Enseig.~Math., Geneva, Switzerland, 1999.

\bibitem[L]{L} J.~Lambek: {\it Lectures on Rings and Modules},
Blaisdell, Waltham, Mass., 1966.

\bibitem[Len]{Len} H.~W.~Lenstra: {\it Complete intersections and Gorenstein
rings}, Proc.~Conf., Elliptic Curves and Modular Forms (Hong Kong, 1993,
J.~Coates and S.~T.~Yau, eds.), 99-109, International Press, 1995.

\bibitem[Le]{Le}  L.~S.~Levy: {\it Modules over Dedekind-like rings},
J.~Algebra {\bf 93}(1985), 1-116.

\bibitem[LR]{LR}  L.~S.~Levy and J.~C.~Robson: {\it Hereditary noetherian 
prime rings III: infinitely generated projective modules}, to appear in
J.~Algebra.

\bibitem[LW]{LW}  L.~S.~Levy and R.~Wiegand: {\it Dedekind-like behavior
of rings with $2$-generated ideals}, J.~Pure Appl.~Algebra 
{\bf 37}(1985), 41-58.

\bibitem[M]{M} E.~Mares: {\it Semiperfect modules}, Math.~Zeit. 
{\bf 82}(1963), 347-360.

\bibitem[Ma]{Ma} J.~Martinet: {\it Modules sur l'alg\`ebre du groupe
quaternionien,} Ann.~Sci.~\'Ecole Norm.~Sup. {\bf 4}(1971), 399-408.

% \bibitem[Mat]{Mat} E.~Matlis: {\it The two-generator problem for ideals},
% Mich. Math. J. {\bf 17}(1970), 157-165.

\bibitem[Me]{Me} P.~Menal: {\it On $\,\pi$-regular rings whose primitive
factor rings are artinian}, J.~Pure Appl.\,Algebra {\bf 20}(1981), 71-81.

\bibitem[MM]{MM} P.~Menal and J.~Moncasi: {\it On regular rings with stable
range} $\,2$, J.~Pure Appl.~Algebra {\bf 24}(1982), 25-40.

\bibitem[Mu]{Mu} M.~P.~Murthy: {\it A survey of obstruction theory for
projective modules of top rank}, this volume, pp.

\bibitem[Mue]{Mue} B.~J.~M\"uller: {\it On semiperfect rings},
Ill. J.~Math. {\bf 14}(1974), 464-467.

\bibitem[NR]{NR} L.~A.~Nazarova and A.~V.~Ro\u{i}ter: {\it Refinement of
a theorem of Bass}, Dokl.~Akad.~Nauk SSSR {\bf 176}(1967), 266-268.
(English Transl.: Soviet Math. {\bf 8}(1967), 1089-1092.)

\bibitem[Ni$_1$]{Ni1} W.~K.~Nicholson: {\it  $I$-rings},  Trans. Amer. 
Math. Soc. {\bf 207}(1975), 361-373.

\bibitem[Ni$_2$]{Ni2} W.~K.~Nicholson: {\it On semiperfect modules},
Canad.~Math.~Bull. {\bf 18}(1975), 77-80.

\bibitem[Ni$_3$]{Ni3} W.~K.~Nicholson: {\it Semiregular modules and rings},
Canad.~J.~Math. {\bf 28}(1976), 1105-1120.

\bibitem[Ni$_4$]{Ni4} W.~K.~Nicholson: {\it Lifting idempotents and exchange 
rings}, Trans. Amer. Math. Soc. {\bf 229}(1977), 269-278.

\bibitem[OS]{OS} U.~Oberst and H.-J.~Schneider: {\it Die Struktur von 
projektiven Moduln}, Invent.~Math {\bf 13}(1971), 295-304.

\bibitem[Os]{Os} B.~Osofsky: {\it A generalization of quasi-Frobenius rings},
J.~Algebra {\bf 4}(1966), 373-387.

\bibitem[Pr]{Pr} M.~Y.~Prest: {\it Duality and pure-semisimple rings},
J.~London Math.~Soc. {\bf 38}(1988), 403-409.

\bibitem[Qu]{Qu}  D.~Quillen: {\it Projective modules over polynomial rings},
Invent. Math. {\bf 36}(1976), 167-171.

\bibitem[Re]{Re} R.~Rentschler: {\it Eine Bemerkung zu Ringen mit 
Minimalbedingung f\"ur Hauptideale}, Arch.~Math. {\bf 17}(1966), 298-301.

\bibitem[Ro]{Ro} A.~V.~Ro\u{i}ter: {\it An analog of Bass' Theorem for
representation modules of noncommutative orders,} Dokl.~Akad.~Nauk SSSR
{\bf 168}(1966), 1261-1264. (English Transl.: Soviet Math.~Dokl. 
{\bf 7}(1966), 830-833.)

\bibitem[Sa$_1$]{Sa1} F.~L.~Sandomierski: {\it On semiperfect and perfect 
rings}, Proc. Amer. Math. Soc. {\bf 21}(1969), 205-207.

\bibitem[Sa$_2$]{Sa2} F.~L.~Sandomierski: {\it Linearly compact modules and
local Morita duality}, in {\it Ring Theory\/} (Proc. Conf., Park City, Utah,
R.~Gordon, ed.), pp.\,333-346, Academic Press, New York, 1972.

\bibitem[Sat]{Sat} H.~Sato: {\it Gorenstein rings with semigroup bases},
J.~Algebra {\bf 88}(1984), 460-488.

\bibitem[S$_1$]{S1} J.-P.~Serre: {\it Faisceaux alg\'ebriques coh\'erents},
Ann. Math. {\bf 61}(1955), 191-278.

\bibitem[S$_2$]{S2} J.-P.~Serre: {\it Modules projectifs et espaces 
fibr\'es \`a fibre vectorielle}, S\'em. Dubreil-Pisot, no.\,23, 
Secretariat Math., Paris, 1957/58.

\bibitem[S$_3$]{S3} J.-P.~Serre: {\it Sur les modules projectifs},
S\'em.~Dubreil-Pisot, no.\,2, Secretariat Math., Paris, 1960/61.

\bibitem[Se$_1$]{Se1} C.~S.~Seshadri: {\it Triviality of vector bundles over 
the affine space $\,K^2$}, Proc. Nat. Acad. Sci. U.S.A. {\bf 44}(1958), 
456-458.

\bibitem[Se$_2$]{Se2} C.~S.~Seshadri: {\it Algebraic vector bundles over
the product of an affine curve and the affine line}, Proc. Amer. Math. Soc.
{\bf 10}(1959), 670-673.

\bibitem[Si$_1$]{Si1} D.~Simson: {\it Functor categories in which every flat 
object is projective}, Bull. Acad. Polon. Sci. {\bf 22}(1974), 375-380.

\bibitem[Si$_2$]{Si2} D.~Simson: {\it Pure semisimple categories and rings
of finite representation type}, J.~Algebra {\bf 48}(1977), 290-296.
(Corrigendum: J.~Algebra {\bf 67}(1980), 254-256.)

\bibitem[Si$_3$]{Si3} D.~Simson: {\it Partial Coxeter functors and right
pure semisimple rings}, J.~Algebra {\bf 71}(1981), 195-218.

\bibitem[Si$_4$]{Si4} D.~Simson: {\it An Artin problem for division ring
extensions and the pure semisimplicity conjecture}, Arch.~Math. {\bf 66}(1996),
114-122.

\bibitem[Sm]{Sm} L.~Small: {\it Reviews in Ring Theory} (as printed
in Math.~Reviews, 1940-79), Amer. Math.~Soc., Providence, R.I., 1981.

%\bibitem[SRS]{SRS} de Smit / Rubin / Schoof
%in: G. Cornell, J.H. Silverman, G. Stevens (eds), Modular
%forms and Fermat's last theorem, Springer, New York, 1997.
% Later work by Diamond gave rise to ring-questions we couldn't solve!

\bibitem[St]{St} J.~T.~Stafford: {\it Stable structure of noncommutative
Noetherian rings}, J.~Algebra {\bf 47}(1977), 244-267.

\bibitem[Su]{Su} A.~Suslin: {\it Projective modules over a polynomial
ring are free}, Dokl. Akad. Nauk SSSR {\bf 229}(1976).  (English
Transl.: Soviet Math. Dokl. {\bf 17}(1976), 160-164.)

\bibitem[Sw$_1$]{Sw1}  R.~Swan: {\it Vector bundles and projective modules},
Trans. Amer. Math. Soc. {\bf 105}(1962), 264-277.

\bibitem[Sw$_2$]{Sw2}  R.~Swan: {\it Projective modules over group rings
and maximal orders}, Annals of Math. {\bf 76}(1962), 55-61.

\bibitem[SE]{SE}  R.~Swan and E.~G.~Evans, Jr.: {\it $\,K$-Theory of
Finite Groups and Orders}, Lecture Notes in Math., Vol.~{\bf 149},
Springer-Verlag, Berlin-Heidelberg-New York, 1970.

\bibitem[Sz]{Sz} F.~Sz\'asz: {\it \"Uber Ringe mit Minimalbedingung f\"ur 
Hauptrechtsideale}, I, Publ. Math. Debrecen {\bf 7}(1960), 54-64; 
II, Acta Math. Acad. Sci. Hungar. {\bf 12}(1961), 417-439;
III, {\it ibid.} {\bf 14}(1963), 447-461.

\bibitem[TW]{TW} R.~Taylor and A.~Wiles: {\it Ring-theoretic properties
of certain Hecke algebras}, Annals of Math.~{\bf 141}(1995), 553-572.

\bibitem[V$_1$]{Va1}  L.~N.~Vaserstein: {\it Stable rank of rings and
dimensionality of topological spaces}, Funct.  Anal. Appl. {\bf 5}(1971),
102-110.

\bibitem[V$_2$]{Va2}  L.~N.~Vaserstein: {\it Bass' first stable range 
condition}, J.~Pure Appl.~Algebra {\bf 34}(1984), 319-330.

\bibitem[Wa$_1$]{Wa1}  R.~Warfield: {\it Exchange rings and decompositions 
of modules}, Math.~Annalen {\bf 199}(1972), 31-36.

\bibitem[Wa$_2$]{Wa2}  R.~Warfield: {\it Cancellation of modules and groups
and stable range of endomorphism rings}, Pac. J. Math {\bf 91}(1980), 457-485.

\bibitem[We]{We}  C.~Weibel: {\it The development of algebraic $K$-theory
before 1980}, this volume, pp.

\bibitem[Wi]{Wi}  R.~Wiegand: {\it Cancellation over commutative rings 
of dimension one and two}, J.~Algebra {\bf 88}(1984), 438-459.

\bibitem[W]{W}  A.~Wiles: {\it Modular elliptic curves and Fermat's
Last Theorem}, Annals of Math. {\bf 142}(1995), 443-551.

\bibitem[Wis]{Wis}  R.~Wisbauer: {\it Foundations of Module and Ring Theory},
Gordon and Breach, Philadelphia, PA, 1991.

\bibitem[Ya]{Ya}  K.~Yamagata: {\it On projective modules with the exchange
property}, Sci. Rep. Tokyo Kyoiku Daigaku, Sect.~A {\bf 12}(1974), 39-48.

\bibitem[Zi$_1$]{Zi1}  B.~Zimmermann-Huisgen: {\it Rings whose right modules
are direct sums of indecomposable modules}, Proc. Amer. Math. Soc. 
{\bf 77}(1979), 191-197.

\bibitem[Zi$_2$]{Zi2}  B.~Zimmermann-Huisgen: {\it On the abundance of
$\,\aleph_1$-separable modules,} in ``Abelian Groups and Noncommutative
Rings'', Contemp.~Math.~{\bf 130}(1992), 167-180, AMS.

\bibitem[ZZ$_1$]{ZZ1}  B.~Zimmermann-Huisgen and W.~Zimmermann: 
{\it Classes of modules with the exchange property}, J.~Algebra 
{\bf 88}(1984), 416-434.

\bibitem[ZZ$_2$]{ZZ2}  B.~Zimmermann-Huisgen and W.~Zimmermann: {\it
On the sparsity of representations of rings of pure global dimension zero},
Trans. Amer. Math. Soc. {\bf 320}(1990), 695-711.

\end{thebibliography}
\end{document}